%% file: Abs_irred_Tate_modules_Arxiv.tex
\documentclass[10pt]{amsart}
\usepackage[margin=1.1in]{geometry}
\usepackage{amsmath,amssymb,amsthm,mathrsfs}
\usepackage[colorlinks=true,linkcolor=blue,citecolor=blue,urlcolor=blue]{hyperref}
\usepackage{tikz}
\usepackage{longtable}
\usepackage{dsfont}
\usepackage{enumitem}
\input{cluster_pictures.tex}

\theoremstyle{plain}
\newtheorem{theorem}{Theorem}[section]
\newtheorem{proposition}[theorem]{Proposition}
\newtheorem{lemma}[theorem]{Lemma}
\newtheorem{corollary}[theorem]{Corollary}
\theoremstyle{definition}
\newtheorem{remark}[theorem]{Remark}
\newtheorem{notation}[theorem]{Notation}
\newtheorem{example}[theorem]{Example}

\newcommand{\Q}{\mathbb{Q}}
\newcommand{\Z}{\mathbb{Z}}
\newcommand{\F}{\mathbb{F}}
\newcommand{\Gal}{\operatorname{Gal}}
\newcommand{\Aut}{\operatorname{Aut}}
\newcommand{\Ind}{\operatorname{Ind}}
\newcommand{\Res}{\operatorname{Res}}
\newcommand{\Jac}{\operatorname{Jac}}
\DeclareMathOperator{\ord}{ord}
\DeclareMathOperator{\Hom}{Hom}
\DeclareMathOperator{\End}{End}

\title[Curves with absolutely irreducible Tate modules over $p$-adic fields]{Curves with absolutely irreducible Tate modules \\ over $p$-adic fields}
\author{Alexandros Konstantinou}
\address{Max Planck Institute for Mathematics, Vivatsgasse 7, 53111 Bonn, Germany}
\email{konstantinou@mpim-bonn.mpg.de}


\begin{document}

\marginpar{}

\begin{abstract}
Over every number field, curves whose Jacobians have absolutely irreducible $\ell$-adic Tate modules exist in every positive genus and can be chosen hyperelliptic. Over $p$-adic fields, however, inertia and Frobenius can exclude entire genera. We classify the genera that occur for $\ell\ne p$ in terms of the residue characteristic and the cardinality of the residue field. For odd $p$, we also classify the genera occurring for hyperelliptic curves. These genera have relative density zero among those occurring for arbitrary curves, yielding Jacobians whose isogeny classes contain no hyperelliptic Jacobian. For odd $p$ and even residue degree, we obtain further restrictions on the genus when the curve itself is required to have potentially good reduction.

Over number fields, we obtain uniform bounds from absolute irreducibility of the Tate module at two finite places of distinct residue characteristics different from $\ell$. The dimension of the abelian variety alone bounds the number of isomorphism classes in its isogeny class and the least degrees of isogenies between them.
\end{abstract}

\maketitle

\setcounter{tocdepth}{1}
\tableofcontents
\section{Introduction}

The paper studies the following question. Given a finite extension $\mathcal K/\Q_p$ and a prime $\ell\ne p$, which positive integers $g$ occur as the genus of a curve over $\mathcal K$ whose Jacobian has an absolutely irreducible $\ell$-adic Tate module?

Over a number field $F$, Faltings' semisimplicity and isogeny theorems show that, for an abelian variety $A/F$, the representation $V_\ell(A)$ is absolutely irreducible if and only if $\End_F(A)=\Z$ \cite{Faltings}. For every $g\ge1$, the hyperelliptic curve $X_g/\Q$ defined by $y^2=x^{2g+1}-x-1$ satisfies $\End_{\overline{\Q}}(\Jac(X_g))=\Z$; see \cite[Thm.~2.1 and Examples~2.2(1)]{Zarhin} for $g\ge2$, while for $g=1$ its non-integral $j$-invariant excludes complex multiplication. Thus, over every number field there are hyperelliptic curves of every positive genus whose Jacobians have absolutely irreducible Tate modules.

Over a finite extension $\mathcal K/\Q_p$, the reduction of an abelian variety imposes further constraints on its Tate module. For instance, if $A/\mathcal K$ has good reduction, inertia acts trivially on $V_\ell(A)$, so the action of $G_{\mathcal K}=\Gal(\overline{\mathcal K}/\mathcal K)$ on $V_\ell(A)$ is determined by a Frobenius lift in $G_{\mathcal K}$. Any eigenvector of this element over $\overline{\Q}_\ell$ spans a $G_{\mathcal K}$-stable subspace, so $V_\ell(A)$ is not absolutely irreducible. Absolute irreducibility therefore forces inertia to act non-trivially, and Frobenius must permute the irreducible constituents of this action transitively. These conditions are restrictive enough to exclude entire genera.

\subsection{Jacobians and local Galois representations}
Our first result classifies the genera of curves over $\mathcal K$ whose Jacobians have absolutely irreducible Tate modules. Here $p$ may be any prime, including 2, and we write $q$ for the cardinality of the residue field of $\mathcal K$.

\begin{theorem}[=Theorems~\ref{thm:classification!} and~\ref{thm:dyadic-jacobian-realisations}]\label{thm:classification}
Let $g\ge1$. There exists a curve $C/\mathcal K$ of genus $g$ such that $V_\ell(\Jac(C))$ is absolutely irreducible and
\begin{itemize}
\item[\textup{(i)}] tamely ramified if and only if $2g=\varphi(m)$ for some integer $m\ge3$ satisfying $(m,p)=1$ and $\ord_m(q)=\varphi(m)$;
\item[\textup{(ii)}] wildly ramified if and only if $(p-1)\mid2g$.
\end{itemize}
Here $\varphi$ denotes Euler's totient function and $\ord_m(q)$ denotes the multiplicative order of $q$ modulo $m$.
\end{theorem}

When $p=2$, the condition in part~\textup{(ii)} is vacuous, so curves whose Jacobians have wildly ramified, absolutely irreducible Tate modules exist over finite extensions of $\Q_2$ in every positive genus. If $\mathcal K$ has even residue degree, then $q$ is a square, so the condition in part~\textup{(i)} cannot hold, since a square cannot generate $(\Z/m\Z)^\times$ when $m\ge3$. Consequently, a curve of genus $g$ whose Jacobian has an absolutely irreducible Tate module exists over such $\mathcal K$ if and only if $(p-1)\mid2g$, and every such Tate module is wildly ramified.

Theorem~\ref{thm:classification} also settles a purely representation-theoretic question. The trace of every element of inertia on the Tate module of an abelian variety is rational \cite[Exp.~IX, Thm.~4.3]{SGA7}. It is therefore natural to ask which dimensions absolutely irreducible $\ell$-adic representations of $G_{\mathcal K}$ satisfying this condition can have. A'Campo gives necessary conditions on these dimensions \cite{ACampo} and, in the tamely ramified case, constructs abstract representations realising them. Remarkably, Theorem \ref{thm:classification} shows that these conditions are sharp for Jacobians in every permitted even dimension. Combining the two gives the following classification:

\begin{corollary}[=Corollary~\ref{cor:abstract}]\label{cor:abstract!}
The dimensions of absolutely irreducible $\ell$-adic representations of $G_{\mathcal K}$ for which the trace of every element of inertia is rational are precisely
\begin{equation*}
\{1\}\cup\{n\ge1:(p-1)\mid n\}\cup
\{\varphi(m):m\ge3,\ (m,p)=1,\ \ord_m(q)=\varphi(m)\}.
\end{equation*}
Every even dimension in this set occurs as the dimension of an absolutely irreducible Tate module $V_\ell(\Jac(C))$ for a curve $C/\mathcal K$. When $\mathcal K$ has even residue degree, the last set in the union is empty.
\end{corollary}

The main result of \cite{ACampo} rules out absolutely irreducible Tate modules for abelian varieties of dimension $7$ when $p>3$. The corollary shows that the hypothesis $p>3$ cannot be dropped, since such Tate modules occur for Jacobians of curves of genus $7$ over every finite extension of $\Q_2$ or $\Q_3$.

\subsection{Hyperelliptic curves and isogeny classes}\label{ssec:hyper}
Over a finite extension $K/\Q_p$ with $p>2$, requiring the curve to be hyperelliptic imposes further restrictions on the genera permitted by Theorem~\ref{thm:classification}, in contrast with the number field case.

\begin{theorem}[=Proposition~\ref{prop:elliptic-classification}, Theorems~\ref{thm:tamehyperelliptic} and~\ref{thm:wildhyperelliptic}, Proposition~\ref{prop:equation-over-K}]\label{thm:hyperelliptic-classification}
Let $g\ge1$. There exists a hyperelliptic curve $C/K$ of genus $g$ such that $V_\ell(\Jac(C))$ is absolutely irreducible and
\begin{itemize}
\item[\textup{(i)}] tamely ramified if and only if one of the following holds:
\begin{itemize}
\item[\textup{(a)}] $g=1$ and either $q\equiv2\pmod3$ or $q\equiv3\pmod4$;
\item[\textup{(b)}] $g\ge2$ and $2g+1$ is prime, distinct from $p$, with $\ord_{2g+1}(q)=2g$;
\end{itemize}
\item[\textup{(ii)}] wildly ramified if and only if $g=(p^a-1)/2$ or $g=p^a-1$ for some integer $a\ge1$.
\end{itemize}
Moreover, every geometrically hyperelliptic\footnote{A curve $C/K$ is \emph{hyperelliptic over $K$} if it admits a morphism of degree two to $\mathbb P^1_K$, or equivalently an equation $y^2=f(x)$ with $f\in K[x]$ square-free, and \emph{geometrically hyperelliptic} if $C_{\overline K}$ is hyperelliptic over $\overline K$. For genus at least two, the two notions differ exactly when the quotient of $C$ by its hyperelliptic involution is a conic without $K$-points.} curve $C/K$ of genus $g\ge2$ such that $V_\ell(\Jac(C))$ is absolutely irreducible is hyperelliptic over $K$.
\end{theorem}

The genera in this classification form a sparse subset of those permitted by general curves in Theorem~\ref{thm:classification}. In the tame case, $2g+1$ must be prime, while in the wild case the genera lie in the two sequences $p^a-1$ and $(p^a-1)/2$, which grow exponentially in $a$. By contrast, every positive multiple of $(p-1)/2$ occurs as the genus of a curve whose Jacobian has an absolutely irreducible Tate module. Since a $K$-isogeny induces an isomorphism of rational Tate modules, these restrictions give an obstruction to a $K$-isogeny class containing a hyperelliptic Jacobian.

To make this comparison precise, for $X\ge1$ let $\mathcal A_K(X)$ be the set of positive integers $g\le X$ for which there exists an abelian variety $A/K$ of dimension $g$ with absolutely irreducible Tate module. Let $\mathcal J_K(X)$ and $\mathcal H_K(X)$ be the subsets obtained by requiring $A$ to be the Jacobian of a curve or of a hyperelliptic curve, respectively. We also write $\mathcal J_K=\bigcup_{X\ge1}\mathcal J_K(X)$ for the set of all genera admitting a curve over $K$ whose Jacobian has an absolutely irreducible Tate module. This is precisely the set described in Theorem~\ref{thm:classification}.

\begin{corollary}[=Corollary~\ref{cor:hyperelliptic-density!}]\label{cor:hyperelliptic-density}
For every $X\ge1$, we have $\mathcal J_K(X)=\mathcal A_K(X)$. Moreover, $\mathcal H_K(X)\subsetneq\mathcal J_K(X)$ for every $X\ge5(p-1)/2$.

Among the genera in $\mathcal J_K$, those admitting hyperelliptic Jacobians with absolutely irreducible Tate modules have relative density zero. More precisely,
\begin{equation*}
\lim_{X\to\infty}\frac{\#\mathcal H_K(X)}{\#\mathcal J_K(X)}=0.
\end{equation*}
Consequently, for almost every $g\in\mathcal J_K$, no Jacobian of dimension $g$ with an absolutely irreducible Tate module is $K$-isogenous to a hyperelliptic Jacobian.
\end{corollary}

Together, Theorems~\ref{thm:classification} and~\ref{thm:hyperelliptic-classification} show that Jacobians realise every dimension available to abelian varieties with absolutely irreducible Tate modules, whereas the genera admitting hyperelliptic Jacobians have relative density zero. The constant $5$ in the bound $X\ge5(p-1)/2$ cannot be replaced by a smaller one uniformly in $K$, since $\mathcal H_{\Q_3}(4)=\mathcal J_{\Q_3}(4)=\{1,2,3,4\}$. The isogeny obstruction already appears in genus $3$, the smallest genus in which a curve can fail to be hyperelliptic. Indeed, the smooth plane quartic $C/\Q_{11}:y^3=x^4-11x$ has a tamely ramified, absolutely irreducible Tate module; see Proposition~\ref{prop:tame-superelliptic-realisation} below. Nonetheless, $\ord_7(11)=3<6$, so Theorem~\ref{thm:hyperelliptic-classification} rules out a hyperelliptic curve of genus $3$ over $\Q_{11}$ with a tamely ramified, absolutely irreducible Tate module. As a result, the $\Q_{11}$-isogeny class of $\Jac(C)$ contains no hyperelliptic Jacobian.

\subsection{Potentially good reduction of curves}

If $V_\ell(\Jac(C))$ is absolutely irreducible, then $\Jac(C)$ has potentially good reduction by Lemma~\ref{lem:arithmetic-consequences}. The curve $C$ itself need not have potentially good reduction, since its stable reduction can be singular even when its Jacobian has good reduction \cite[\S\S1.2--1.3]{MoonenOort}. We show that some Jacobians are not $K$-isogenous to the Jacobian of any curve with potentially good reduction.

In the tame case, requiring the curve to have potentially good reduction imposes no additional restriction on the genus. The superelliptic constructions of \S\ref{ssec:tame-abelian-realisations} have potentially good reduction and realise every genus permitted by Theorem~\ref{thm:classification}\textup{(i)}. When $K$ has even residue degree, however, $q$ is a square, so every absolutely irreducible Jacobian Tate module is wildly ramified. In this setting, the additional requirement that the curve have potentially good reduction excludes an interval of genera.

\begin{theorem}[=Theorem~\ref{thm:gap!}]\label{cor:gap}
Suppose that $K$ has even residue degree. Let $C/K$ be a curve of genus $g\ge1$ with potentially good reduction. If $V_\ell(\Jac(C))$ is absolutely irreducible, then
\begin{equation*}
g=\frac{p-1}{2}
\qquad\text{or}\qquad
g\ge\frac{p(p-1)}{2}.
\end{equation*}
\end{theorem}

Both endpoints of this gap are attained for suitable $K$; see Proposition~\ref{prop:wild-single-cluster} for the lower endpoint and Lemma~\ref{lem:gap-sharpness} for the upper endpoint.

For every finite extension $K/\Q_p$ of even residue degree and every $2\le r\le p-1$, Theorem~\ref{thm:classification}\textup{(ii)} provides a curve $C/K$ of genus $r(p-1)/2$ whose Jacobian has an absolutely irreducible Tate module. These genera lie strictly inside the excluded interval. Since $K$-isogenies preserve the dimension and the rational Tate module, the theorem rules out potentially good reduction for every curve whose Jacobian is $K$-isogenous to $\Jac(C)$, as stated in Corollary~\ref{cor:no-potential-good-isogeny!}. Thus the distinction between potentially good reduction of a curve and of its Jacobian can persist throughout an entire $K$-isogeny class.

\subsection{Uniform isogeny bounds over number fields}

We also obtain finiteness results for isogenies over number fields. As part of Faltings' work on the Mordell conjecture, he proved that every $M$-isogeny class of abelian varieties over a number field $M$ contains only finitely many $M$-isomorphism classes \cite{Faltings}. Masser and W\"ustholz \cite{MasserWustholz} bounded the least degree of an $M$-isogeny in terms of the dimension, the Faltings height and $[M:\Q]$, with explicit bounds given in \cite[Thm.~1.4]{GaudronRemond}.

For elliptic curves over $\Q$, Mazur classified the prime degrees of rational isogenies \cite{Mazur}, and Kenku proved that each $\Q$-isogeny class contains at most eight $\Q$-isomorphism classes \cite{Kenku}. The following theorem gives bounds in every dimension under a condition at two finite places. These bounds depend only on the dimension and are independent of the height and the number field.

For an abelian variety $A/M$ over a number field, we say that $V_\ell(A)$ is \emph{locally irreducible at $v$}, where $v\nmid\ell$ is a finite place of $M$, if its restriction to $G_{M_v}$ is absolutely irreducible. If this holds for one prime $\ell$ different from the residue characteristic of $v$, then $A$ has potentially good reduction at $v$, and \cite[Cor.~1.6]{DDWeil} implies that it holds for every such prime.

\begin{theorem}[=Theorem~\ref{thm:uniform-isogeny-bound!}]\label{thm:uniform-isogeny-bound}
Let $A/M$ be an abelian variety of dimension $g\ge1$ over a number field. Suppose that, for some prime $\ell$, the representation $V_\ell(A)$ is locally irreducible at two finite places $v,w\nmid\ell$ of distinct residue characteristics.

Then there exists an integer $B_g\ge1$, depending only on $g$, such that every abelian variety $B/M$ that is $M$-isogenous to $A$ admits an $M$-isogeny $A\to B$ of degree at most $B_g$. Moreover, the $M$-isogeny class of $A$ contains at most $B_g$ $M$-isomorphism classes.
\end{theorem}

For Jacobians, this also bounds the number of curves of a fixed genus whose Jacobians belong to a fixed isogeny class. More precisely, if $C/M$ has genus $g\ge2$ and $\Jac(C)$ satisfies the hypotheses above, then the number of $M$-isomorphism classes of curves $D/M$ with $\Jac(D)$ $M$-isogenous to $\Jac(C)$ is bounded solely in terms of $g$; see Corollary~\ref{cor:uniform-curves-isogeny-class} below.

\subsection{Methods and organisation of the paper}\label{ssec:methods}

To prove the equivalences of Theorems~\ref{thm:classification} and~\ref{thm:hyperelliptic-classification}, we derive necessary conditions and construct curves realising every permitted genus. Both directions rest on the action of inertia and Frobenius on the Tate module, studied in \S\ref{sec:local-representations}. By Lemma~\ref{lem:arithmetic-consequences}, an abelian variety with an absolutely irreducible Tate module has potentially good reduction, so inertia acts through a finite quotient.

In \S\ref{ssec:inertia-automorphisms}, we recall how reduction allows us to describe this action geometrically. When a curve $C/\mathcal K$ itself acquires good reduction over a finite Galois extension $F/\mathcal K$, reduction identifies $V_\ell(\Jac(C))$ with $V_\ell(\Jac(\overline C))$ as $I_{\mathcal K}$-representations, where inertia acts on the geometric special fibre $\overline C$ by automorphisms (Proposition~\ref{prop:inertia-automorphisms}). We then use two facts to study $V=V_\ell(\Jac(C))\otimes_{\Q_\ell}\overline{\Q}_\ell$ as an $I_{\mathcal K}$-representation.
\begin{enumerate}[label=\textup{(Fact~\arabic*)},ref=\arabic*,leftmargin=*,labelindent=\parindent]
\item\label{fact:traces} The trace of every element of $I_{\mathcal K}$ on $V$ is rational.
\item\label{fact:invariants} For every subgroup $H\subseteq I(F/\mathcal K)$, we have $\dim V^H=2g(\overline C/H)$.
\end{enumerate}
It remains to determine how Frobenius permutes the irreducible constituents of $V|_{I_{\mathcal K}}$ (Lemma~\ref{lem:frobenius-transitivity}).

\subsubsection*{Jacobians}
Our constructions start with two elliptic curves over $\Q_3$.
\begin{example}\label{ex:tame-elliptic}
Consider $E/\Q_3:y^2=x^3-3x$. Choose $t$ with $t^4=3$ and put $F=\Q_3(t,\mu_4)$. This is a tame Galois extension, and the substitution $x=t^2u$, $y=t^3w$ gives the good-reduction model $w^2=u^3-u$. The inertia group $I(F/\Q_3)$ is cyclic of order four, and a generator induces the automorphism $\gamma:(u,w)\mapsto(-u,iw)$ on the geometric special fibre $\overline E:w^2=u^3-u$, where $i^2=-1$.

Fix $\ell\ne3$ and put $V=V_\ell(E)\otimes_{\Q_\ell}\overline{\Q}_\ell$ and $I=\langle\gamma\rangle$. Since $\gamma^2$ is the hyperelliptic involution, Fact~\ref{fact:invariants} gives $V^{\langle\gamma^2\rangle}=0$. Thus only the two faithful characters $\chi_1,\chi_2$ of $I$ can occur in $V|_I$. By Fact~\ref{fact:traces}, their multiplicities are equal, and since $\dim V=2$, we obtain $V|_I\cong\chi_1\oplus\chi_2$. An arithmetic Frobenius lift $\Phi$ satisfies $\Phi\gamma\Phi^{-1}=\gamma^{-1}$ and therefore exchanges the two characters. Thus $V_\ell(E)$ is tamely ramified and absolutely irreducible by Lemma~\ref{lem:frobenius-transitivity}.
\end{example}

\begin{example}\label{ex:wild-elliptic}
Let $b\in\Z_3^\times$ and consider $E_b/\Q_3:y^2=x^3-3bx-3$. Choose a root $\alpha$ of $x^3-3bx-3$ and $t$ with $t^4=3b$. The substitution $x=\alpha+t^2u$, $y=t^3w$ gives good reduction over the Galois extension $F=\Q_3(\alpha,t,\mu_4)$, with geometric special fibre $\overline E_b:w^2=u^3-u$. The inertia group $I(F/\Q_3)\cong C_3\rtimes C_4$ acts faithfully on $\overline E_b$, with generators inducing $\tau:(u,w)\mapsto(u+1,w)$ and $\gamma:(u,w)\mapsto(-u,iw)$, where $i^2=-1$ and $\gamma\tau\gamma^{-1}=\tau^{-1}$.

Fix $\ell\ne3$ and put $V=V_\ell(E_b)\otimes_{\Q_\ell}\overline{\Q}_\ell$ and $P=\langle\tau\rangle$. The quotient $\overline E_b/P$ has genus zero with parameter $w$, so Fact~\ref{fact:invariants} gives $V^P=0$. Thus only the two non-trivial characters $\psi_1,\psi_2$ of $P$ can occur in $V|_P$. By Fact~\ref{fact:traces}, their multiplicities are equal, and since $\dim V=2$, we obtain $V|_P\cong\psi_1\oplus\psi_2$. The element $\gamma$ exchanges the two characters, so $V$ is already irreducible under inertia. Hence $V_\ell(E_b)$ is wildly ramified and absolutely irreducible.
\end{example}

We extend these examples to higher genus in \S\S\ref{sec:hyperelliptic-realisations}--\ref{sec:abelian-realisations}. Generalising Example~\ref{ex:tame-elliptic}, the superelliptic curves $Y_{\nu,a}:y^\nu=x^{\nu^{a-1}+1}-\pi x$ realise every genus permitted by Theorem~\ref{thm:classification}\textup{(i)}. Example~\ref{ex:wild-elliptic} generalises to the hyperelliptic curves $X_b:y^2=x^p-\pi bx-\pi$, which realise the least wild genus $(p-1)/2$. This case lies in the family studied in \cite[Prop.~2.2 and Thm.~2.3]{Coppola}; the constructions below also treat the higher powers $p^a$ and the two-cluster case. To obtain every genus permitted by Theorem~\ref{thm:classification}\textup{(ii)}, we apply Weil restriction from unramified extensions and then use a smoothing argument to realise the resulting Tate modules by Jacobians. The constructions over finite extensions of $\Q_2$ follow the same strategy (\S\ref{ssec:dyadic-realisations}).

\subsubsection*{Hyperelliptic curves}
For $K/\Q_p$ with $p>2$, we use the cluster pictures of \cite{DDMM}, recalled in \S\ref{ssec:hyperelliptic-clusters}, to obtain the necessary restrictions in \S\ref{sec:hyperelliptic-irreducibility}. For $g\ge2$, Corollary~\ref{cor:cluster-configurations} shows that absolute irreducibility forces a suitable equation to have either a single principal cluster or two Galois-conjugate principal clusters. The action of inertia on their roots then restricts the genus. These configurations also guide the constructions in \S\ref{sec:hyperelliptic-realisations}, which realise every genus in Theorem~\ref{thm:hyperelliptic-classification}.

\begin{center}
\renewcommand{\arraystretch}{1.3}
\begin{tabular}{@{}lccc@{}}
\hline
Curve & Genus & Ramification & Principal clusters\\
\hline
$X_{a,b,\eta}:y^2=x^{p^a}-\pi bx-\pi\eta$
&$(p^a-1)/2$ & Wild & one\\
$D_a:y^2=f_\delta(x)f_{-\delta}(x)$
&$p^a-1$ & Wild & two\\
$H_r:y^2=x^r-\pi$
&$(r-1)/2$ & Tame & one\\
$E:y^2=x^3-\pi x$
&$1$ & Tame & one\\
\hline
\end{tabular}
\end{center}

Here $a\ge1$ and $b,\eta\in\mathcal O_K^\times$. The polynomial $f_{\pm\delta}$ is the one defining $X_{a,\pm\delta,1}$ over $K(\delta)$, translated to the centre $\pm\delta$, where $\delta^2\in\mathcal O_K^\times$ has non-square reduction. The two wild families are treated in \S\S\ref{ssec:wild-single-cluster}--\ref{ssec:wild-two-clusters}, respectively. In the tame rows, $r\ne p$ is a prime with $r\ge5$ and $\ord_r(q)=r-1$. In genus one, $E$ is used when $q\equiv3\pmod4$, while the case $q\equiv2\pmod3$ uses the elliptic curve $y^2-\pi y=x^3$; {see Proposition~\ref{prop:elliptic-classification}.} The family $H_r$ is treated in \S\ref{sec:tame-hyperelliptic-realisations}.

\medskip

In \S\ref{sec:genus-gap}, we consider curves with potentially good reduction and prove the genus gap of Theorem~\ref{thm:gap!}. Finally, in \S\ref{sec:global}, we obtain bounds over number fields for isogeny degrees and the number of isomorphism classes in an isogeny class. These depend only on the dimension and assume local irreducibility at two finite places of distinct residue characteristics.

\subsection*{Acknowledgements}

This work began while I was supported by the scientific programme ``Enhancing the Research Capacity in Mathematical Sciences (PIKOM)'', No.~DO1-67/05.05.2022, of the Ministry of Education and Science of Bulgaria. I am grateful to the Max Planck Institute for Mathematics in Bonn for its hospitality, financial support and excellent working conditions. I thank the MathOverflow community, and in particular the user NMarkovBlanket, for the argument used in the proofs of Corollary~\ref{cor:wild-jacobian-realisations} and Theorem~\ref{thm:dyadic-jacobian-realisations}; see \cite{MO515765}. In addition, I used ChatGPT (Sol 6.1) to obtain comments on a first complete draft of this paper. In particular, it helped develop the argument presented on MathOverflow into a complete proof of Lemma~\ref{lem:compact-type-smoothing}. Its comments also suggested the result now stated in Proposition~\ref{prop:equation-over-K}, which I subsequently formalised.

\subsection{Notation}

We write $\mathcal K$ for a finite extension of $\Q_p$, with no restriction on the prime $p$. We reserve $K$ for finite extensions of $\Q_p$ with $p>2$, and $L$ for finite extensions of $\Q_2$. Over a local field of residue characteristic $p$, we always assume that $\ell\ne p$. In the global applications, $M$ denotes a number field and $M_v$ its completion at a finite place $v$.

We use the following notation.

\begingroup
\renewcommand{\arraystretch}{1.2}
\begin{longtable}{@{}lp{0.66\textwidth}@{}}
\hline
Notation & Meaning\\
\hline
\multicolumn{2}{@{}l}{\textit{Local fields}}\\
$\mathcal O_{\mathcal K},\mathfrak m_{\mathcal K}$
&Ring of integers of $\mathcal K$ and its maximal ideal\\
$k,q$
&Residue field of $\mathcal K$ and its cardinality\\
$\pi,v$
&A uniformiser of $\mathcal K$ and the valuation normalised by $v(\pi)=1$\\
$\overline{\mathcal K},\overline k$
&A fixed algebraic closure of $\mathcal K$ and its residue field\\
\hline
\multicolumn{2}{@{}l}{\textit{Galois groups and representations}}\\
$G_M$
&Absolute Galois group $\Gal(\overline M/M)$ of a field $M$\\
$I_M,P_M$
&Inertia and wild inertia subgroups of $G_M$ when $M$ is a local field\\
$I(F/\mathcal K),P(F/\mathcal K)$
&Inertia and wild inertia subgroups of a finite Galois extension $F/\mathcal K$\\
$G,I,P$
&Images of $G_{\mathcal K},I_{\mathcal K},P_{\mathcal K}$ under $\rho$\\
$\operatorname{Irr}(H)$
&Set of isomorphism classes of irreducible $\overline{\Q}_\ell$-representations of a finite group $H$\\
$\widehat H$
&Character group $\Hom(H,\overline{\Q}_\ell^\times)$ of a finite abelian group $H$\\
$\mathds{1},V^H$
&Trivial representation and subspace of $H$-invariants in $V$\\
\hline
\multicolumn{2}{@{}l}{\textit{Curves and abelian varieties}}\\
$g(C),\Jac(C)$
&Genus and Jacobian of a curve $C$\\
$T_\ell(A),V_\ell(A)$
&Integral and rational $\ell$-adic Tate modules of an abelian variety $A$\\
$\Res_{K'/K}$
&Restriction of scalars along a finite extension $K'/K$\\
\hline
\multicolumn{2}{@{}l}{\textit{Groups and arithmetic}}\\
$C_n,\mu_n$
&Cyclic group of order $n$ and group of $n$th roots of unity\\
$C_J(H)$
&Centraliser of a subgroup $H$ in $J$\\
$\F_{p^a}$
&Finite field of cardinality $p^a$\\
$\varphi(n),\ord_n(q)$
&Euler's totient function and the multiplicative order of $q$ modulo $n$\\
$\nu$
&A rational prime distinct from $p$\\
\hline
\end{longtable}
\endgroup

\subsubsection*{Local fields}
We extend the chosen valuation on each local field to its fixed algebraic closure. An arithmetic Frobenius lift is an element $\Phi\in G_{\mathcal K}$ that induces $x\mapsto x^q$ on $\overline k$. We use $\ord_n(q)$ only when $(n,q)=1$, with the convention $\ord_1(q)=1$.

\subsubsection*{Curves and abelian varieties.}
Unless otherwise stated, curves are smooth, projective and geometrically connected. If $C$ acquires good reduction over a finite extension $F/\mathcal K$, then $\overline C$ denotes the geometric special fibre of its smooth proper model over $\mathcal O_F$. For an abelian variety $A$, we put
\begin{equation*}
V=V_\ell(A)\otimes_{\Q_\ell}\overline{\Q}_\ell,
\end{equation*}
with $A=\Jac(C)$ when studying a curve $C$.

\subsubsection*{Representations and characters.}
All $\ell$-adic representations are finite-dimensional and continuous. We use $W$ for a general $\Q_\ell$-representation and $\rho$ for the representation under consideration. We say that $\rho$ has rational traces on inertia if $\operatorname{tr}(\rho(\sigma))\in\Q$ for every $\sigma\in I_{\mathcal K}$.

When considering characters and irreducible constituents, we work over $\overline{\Q}_\ell$. For a finite group $H$, an irreducible representation $\chi\in\operatorname{Irr}(H)$ and a $\overline{\Q}_\ell$-representation $U$ of $H$, we write $\langle\chi,U\rangle_H$ for the multiplicity of $\chi$ in $U$.

For a one-dimensional character $\chi$ of a normal subgroup $H\subseteq G$, we define $(\sigma\cdot\chi)(h)=\chi(\sigma^{-1}h\sigma)$ for $\sigma\in G$ and $h\in H$. We write $\chi^n$ for integer powers of a character and $\chi^{\oplus n}$ for a direct sum of $n$ copies.

\section{Local Galois representations}\label{sec:local-representations}

We study the restriction to inertia of an absolutely irreducible $G_{\mathcal K}$-representation and recall the dimension restrictions when its traces on inertia are rational. Finally, we study what absolute irreducibility of the Tate module implies for the endomorphisms and reduction of an abelian variety.

Throughout this section, let $\mathcal K/\Q_p$ be a finite extension, where $p$ is any prime, and let $\ell\ne p$ be prime. Write $q$ for the cardinality of the residue field of $\mathcal K$.

\subsection{Rational traces on inertia and dimension restrictions}

For a non-zero finite-dimensional continuous representation $\rho:G_{\mathcal K}\to\mathrm{GL}(W)$ over $\Q_\ell$, we consider the following conditions:
\begin{itemize}
\item[(R1)] $W$ is absolutely irreducible;
\item[(R2)] $\operatorname{tr}(\rho(\sigma))\in\Q$ for every $\sigma\in I_{\mathcal K}$.
\end{itemize}

Write $G,I,P$ for the images of $G_{\mathcal K},I_{\mathcal K},P_{\mathcal K}$ under $\rho$. By the monodromy argument in \cite[proof of Prop.~3]{ACampo}, condition \textup{(R1)} implies that $I$ is finite. This argument does not require \textup{(R2)} or tame ramification. It follows from Maschke's theorem that the restriction of $W\otimes_{\Q_\ell}\overline{\Q}_\ell$ to $I$ is semisimple.

We write $\operatorname{Irr}(I)$ for the set of isomorphism classes of irreducible $\overline{\Q}_\ell$-representations of $I$. An arithmetic Frobenius lift $\Phi$ acts on $\operatorname{Irr}(I)$ by
\begin{equation*}
(\Phi\cdot\chi)(h)=\chi(\Phi^{-1}h\Phi),
\qquad h\in I.
\end{equation*}
Thus $\Phi$ sends the $\chi$-isotypic component of $W\otimes_{\Q_\ell}\overline{\Q}_\ell$ to the $(\Phi\cdot\chi)$-isotypic component.

Since inertia and $\Phi$ generate a dense subgroup of $G_{\mathcal K}$, continuity implies that a subspace stable under both is stable under $G_{\mathcal K}$. The following lemma shows that absolute irreducibility forces each irreducible constituent of the restriction to inertia to occur with multiplicity one, and Frobenius to permute these constituents transitively.

\begin{lemma}\label{lem:multiplicity-free}
Suppose that $W$ is absolutely irreducible. Then there exists a subset $\Omega\subseteq\operatorname{Irr}(I)$ such that
\begin{equation*}
(W\otimes_{\Q_\ell}\overline{\Q}_\ell)|_I\cong\bigoplus_{\chi\in\Omega}\chi,
\end{equation*}
and $\Phi$ acts transitively on $\Omega$.
\end{lemma}

\begin{proof}
Let $\Phi$ denote the action of an arithmetic Frobenius lift on $W\otimes_{\Q_\ell}\overline{\Q}_\ell$, and put
\begin{equation*}
\Omega=\{\chi\in\operatorname{Irr}(I):\langle\chi,W\otimes_{\Q_\ell}\overline{\Q}_\ell\rangle_I\ne0\},
\qquad
r=|\Omega|.
\end{equation*}
For each $\chi\in\Omega$, let $W_\chi$ be the sum of all $I$-submodules of $W\otimes_{\Q_\ell}\overline{\Q}_\ell$ which are isomorphic to $\chi$. Thus $W_\chi$ is the $\chi$-isotypic component, and $W_\chi\cong\chi^{\oplus n_\chi}$ for some $n_\chi\ge1$.

The element $\Phi$ permutes these isotypic components transitively. Otherwise, the direct sum of the components in one orbit would be a proper non-zero $G_{\mathcal K}$-stable subspace, contrary to absolute irreducibility. Fix $\chi\in\Omega$. Then
\begin{equation*}
W\otimes_{\Q_\ell}\overline{\Q}_\ell=\bigoplus_{j=0}^{r-1}\Phi^jW_\chi,
\qquad
\Phi^rW_\chi=W_\chi.
\end{equation*}
For a representation $\chi\in\Omega$, write $\rho_\chi$ for its representation of $I$. Since conjugation by $\Phi^r$ preserves the isomorphism class of $\chi$, there exists $B\in\operatorname{GL}(\chi)$ satisfying $B\rho_\chi(h)B^{-1}=\rho_\chi(\Phi^rh\Phi^{-r})$ for every $h\in I$.

Consider the linear map
\begin{equation*}
T:\Hom_I(\chi,W_\chi)\to\Hom_I(\chi,W_\chi),
\qquad
T(f)=\Phi^r\circ f\circ B^{-1}.
\end{equation*}
The above identity for $B$ gives $T(f)\rho_\chi(h)=hT(f)$ for every $h\in I$, so $T$ is well defined. Its inverse sends $f$ to $\Phi^{-r}\circ f\circ B$. Since the coefficient field is algebraically closed, choose a non-zero eigenvector $f_\chi:\chi\to W_\chi$ satisfying $T(f_\chi)=\lambda f_\chi$. The eigenvalue $\lambda$ is non-zero because $T$ is invertible. The identity $\Phi^rf_\chi(v)=\lambda f_\chi(Bv)$ for every $v\in\chi$ shows that $\Phi^r$ preserves $f_\chi(\chi)$. Furthermore, $f_\chi$ is injective because it is a non-zero $I$-equivariant map from the irreducible representation $\chi$.

It follows that $W'=\bigoplus_{j=0}^{r-1}\Phi^jf_\chi(\chi)$ is a non-zero subspace stable under $I$ and $\Phi$, and hence under $G_{\mathcal K}$. Since $W$ is absolutely irreducible, $W'=W\otimes_{\Q_\ell}\overline{\Q}_\ell$. Comparing the $\chi$-isotypic components gives $f_\chi(\chi)=W_\chi$. Thus $W_\chi\cong\chi$ and $n_\chi=1$. Transitivity gives the same conclusion for every constituent.
\end{proof}

Conversely, when inertia acts through a finite group, multiplicity one and transitivity under Frobenius imply absolute irreducibility.

\begin{lemma}\label{lem:frobenius-transitivity}
Let $W$ be a non-zero finite-dimensional continuous $\Q_\ell$-representation of $G_{\mathcal K}$ on which inertia acts through a finite group $I$. Suppose that
\begin{equation*}
(W\otimes_{\Q_\ell}\overline{\Q}_\ell)|_I\cong\bigoplus_{\chi\in\Omega}\chi
\end{equation*}
for a subset $\Omega\subseteq\operatorname{Irr}(I)$. Then $W$ is absolutely irreducible if and only if an arithmetic Frobenius lift $\Phi$ permutes these irreducible summands transitively.
\end{lemma}

\begin{proof}
Since the restriction to $I$ is semisimple and its irreducible summands are pairwise non-isomorphic, every $I$-stable subspace of $W\otimes_{\Q_\ell}\overline{\Q}_\ell$ is a direct sum of some of these summands. Such a subspace is stable under $\Phi$ precisely when the corresponding subset of $\Omega$ is $\Phi$-stable. Since inertia and $\Phi$ generate a dense subgroup of $G_{\mathcal K}$, it follows that $W$ is absolutely irreducible if and only if $\Phi$ acts transitively \hbox{on $\Omega$}.
\end{proof}

\begin{remark}\label{rem:inertia-transitivity}
The same argument applies to a normal subgroup $P$ of $I$. If $(W\otimes_{\Q_\ell}\overline{\Q}_\ell)|_P$ is a direct sum of pairwise non-isomorphic irreducible representations, then $W|_I$ is absolutely irreducible if and only if $I$ permutes these summands transitively.
\end{remark}

Lemma~\ref{lem:multiplicity-free} uses only condition \textup{(R1)}. The following example shows that this condition does not imply \textup{(R2)}.

\begin{example}
Since $\Q_7$ contains the cube roots of unity, the Eisenstein polynomial $x^3-7$ defines a totally ramified cyclic extension $N=\Q_7(\sqrt[3]{7})$ of degree $3$. Choose a generator $\sigma$ of $\Gal(N/\Q_7)$ and a primitive cube root of unity $\zeta_3\in\Q_{13}$. Consider the character
\begin{equation*}
\chi:G_{\Q_7}\twoheadrightarrow\Gal(N/\Q_7)\longrightarrow\Q_{13}^\times,
\qquad
\sigma\mapsto\zeta_3.
\end{equation*}
This one-dimensional representation satisfies \textup{(R1)}. Since $N/\Q_7$ is totally ramified, some $h\in I_{\Q_7}$ maps to $\sigma$, giving $\chi(h)=\zeta_3\notin\Q$. Thus \textup{(R2)} fails.
\end{example}

Suppose that \textup{(R1)} holds. Both $I$ and $P$ are normal in $G$, the group $P$ is the unique Sylow $p$-subgroup of $I$, and $I/P$ is cyclic of order prime to $p$. We say that $W$ is \emph{tamely ramified} if $P=1$, and \emph{wildly ramified} otherwise.

\begin{lemma}[\cite{ACampo}]\label{lem:tame}
Let $W$ satisfy \textup{(R1)--(R2)}. Then the following hold.
\begin{itemize}
\item[\textup{(i)}] If $P=1$, then $I$ is cyclic of order $m$ prime to $p$, and $\dim W=\varphi(m)=\ord_m(q)$, where \hbox{$\ord_1(q)=1$.}
\item[\textup{(ii)}] If $P\ne1$, then $(p-1)\mid\dim W$.
\end{itemize}
In particular, if $q$ is a square, then either $\dim W=1$, or $W$ is wildly ramified and $(p-1)\mid\dim W$.
\end{lemma}

\begin{proof}
Parts \textup{(i)} and \textup{(ii)} follow from \cite[proof of Prop.~3 and Prop.~4]{ACampo}, respectively. If $q$ is a square, then $\ord_m(q)=\varphi(m)$ forces $m\in\{1,2\}$, since a square cannot generate $(\Z/m\Z)^\times$ when $m\ge3$. The final assertion follows.
\end{proof}

\subsection{$\ell$-adic Tate modules of abelian varieties}\label{ssec:tate-modules}

The following lemma establishes \textup{(R2)} for Tate modules of abelian varieties and records the consequences of absolute irreducibility for endomorphisms and reduction.

\begin{lemma}\label{lem:arithmetic-consequences}
Let $A/\mathcal K$ be an abelian variety of positive dimension and suppose that $V_\ell(A)$ is absolutely irreducible. Then the following hold.
\begin{itemize}
\item[\textup{(i)}] $A$ is simple over $\mathcal K$ and $\End_{\mathcal K}(A)=\Z$;
\item[\textup{(ii)}] $A$ has potentially good reduction and purely additive reduction over $\mathcal K$;
\item[\textup{(iii)}] $V_\ell(A)$ satisfies \textup{(R2)}.
\end{itemize}
\end{lemma}

\begin{proof}
If $B\subseteq A$ is an abelian subvariety defined over $\mathcal K$, then $V_\ell(B)$ is a $G_{\mathcal K}$-stable subspace of $V_\ell(A)$, so $B=0$ or $B=A$. Moreover, absolute irreducibility gives $\End_{G_{\mathcal K}}(V_\ell(A))=\Q_\ell.$
Since $\End_{\mathcal K}(A)\otimes_{\Z}\Q_\ell$ embeds in $\End_{G_{\mathcal K}}(V_\ell(A))$, we obtain $\End_{\mathcal K}(A)\otimes_{\Z}\Q=\Q$. The ring $\End_{\mathcal K}(A)$ is an order in $\Q$, so it equals $\Z$. This proves \textup{(i)}.

As observed above, absolute irreducibility implies that inertia acts through a finite quotient. By the N\'eron--Ogg--Shafarevich criterion, $A$ has potentially good reduction \cite[Thms.~1 and~2(i)]{SerreTate}. The subspace $V_\ell(A)^{I_{\mathcal K}}$ is $G_{\mathcal K}$-stable, so it is $0$ or $V_\ell(A)$. If $V_\ell(A)^{I_{\mathcal K}}=V_\ell(A)$, then $V_\ell(A)$ would be an absolutely irreducible representation of the abelian group $G_{\mathcal K}/I_{\mathcal K}$. All such representations are one-dimensional, whereas $\dim V_\ell(A)=2\dim A\ge2$. Thus $V_\ell(A)^{I_{\mathcal K}}=0$.

If $a$ and $t$ denote the abelian and toric ranks of the identity component of the special fibre of the N\'eron model, then
\begin{equation*}
\dim V_\ell(A)^{I_{\mathcal K}}=2a+t
\end{equation*}
by \cite[\S1, Lems.~1 and~2]{SerreTate}. Hence $a=t=0$, which means that the reduction is purely additive. This proves \textup{(ii)}.

Since $A$ has potentially good reduction, part~\textup{(iii)} follows from \cite[Thm.~2(ii)]{SerreTate}.
\end{proof}

In fact, $A$ acquires supersingular reduction. More precisely, Frobenius acts as an integer scalar over every finite Galois extension over which $A$ has good reduction.

\begin{proposition}\label{prop:square-residue-good-reduction}
Let $A/\mathcal K$ be an abelian variety with absolutely irreducible Tate module $V_\ell(A)$, and let $F/\mathcal K$ be a finite Galois extension over which $A$ has good reduction $\widetilde A/k_F$. Then any arithmetic Frobenius lift in $G_F$ acts on $V_\ell(A)$ as a scalar $\lambda\in\Z$ satisfying $\lambda^2=|k_F|$. In particular, $|k_F|$ is a square and $\widetilde A$ is supersingular.
\end{proposition}

\begin{proof}
Write $\rho$ for the representation on $V_\ell(A)$, and let $\Phi_F\in G_F$ be an arithmetic Frobenius lift. Since $F/\mathcal K$ is Galois, $G_F$ is normal in $G_{\mathcal K}$. For every $\sigma\in G_{\mathcal K}$, the element $\sigma\Phi_F\sigma^{-1}\Phi_F^{-1}$ belongs to $G_F$. It also belongs to $I_{\mathcal K}$, since $G_{\mathcal K}/I_{\mathcal K}$ is abelian, and hence lies in $G_F\cap I_{\mathcal K}=I_F$. Good reduction and the N\'eron--Ogg--Shafarevich criterion imply that $\rho(I_F)=1$, so $\rho(\Phi_F)$ commutes with $\rho(G_{\mathcal K})$.

Schur's lemma over $\overline{\Q}_\ell$ implies that $\rho(\Phi_F)$ is a scalar. Since this endomorphism is defined over $\Q_\ell$, the scalar $\lambda$ belongs to $\Q_\ell$.

Under reduction, $\rho(\Phi_F)$ identifies with the action of the Frobenius endomorphism of $\widetilde A$, whose characteristic polynomial lies in $\Z[T]$; see \cite[Thm.~3]{SerreTate}. Since $\rho(\Phi_F)$ acts by multiplication by $\lambda$, its characteristic polynomial is $(T-\lambda)^{2\dim A}\in\Z[T]$. Its trace is $2\dim(A)\lambda\in\Z$, so $\lambda\in\Q$. Since $\lambda$ is also an algebraic integer, $\lambda\in\Z$. The Weil bounds give $\lambda^2=|k_F|$.

Every Frobenius eigenvalue has $p$-adic valuation equal to half that of $|k_F|$, so all Newton slopes are $1/2$. Thus $\widetilde A$ is supersingular; see \cite[Def.~4.1 and Thm.~4.2]{Oort}.
\end{proof}

\section{Reduction of curves and Galois actions}\label{sec:reduction-galois}

We recall how reduction describes the Galois action on the Tate module of a Jacobian. We first consider curves with potentially good reduction and describe the inertia action through automorphisms of the special fibre, following \cite{DDM}. We then introduce cluster pictures for hyperelliptic curves and their description of the Tate module, following \cite{DDMM}. For a curve $C/\mathcal K$, write $V=V_\ell(\Jac(C))\otimes_{\Q_\ell}\overline{\Q}_\ell$.

We first record a lemma on finite group actions on curves over $\overline k$ which will be used repeatedly in the sequel.

\begin{lemma}\label{lem:quotient-invariants}
Let $X/\overline k$ be a curve and let $H\subseteq\Aut_{\overline k}(X)$ be finite. Then:
\begin{itemize}
\item[\textup{(i)}] The $H$-representation $V_\ell(\Jac(X))$ has rational character.
\item[\textup{(ii)}] The dimension of its $H$-invariant subspace is $\dim V_\ell(\Jac(X))^H=2g(X/H)$.
\end{itemize}
\end{lemma}

\begin{proof}
The character is rational, since endomorphisms of $\Jac(X)$ have characteristic polynomials in $\Z[T]$ \cite[\S19, Thm.~4]{Mumford}. The argument for the first assertion of \cite[Lem.~4.28]{KonstantinouMorgan} applies unchanged in characteristic $p$ and shows that pullback along $X\to X/H$ gives an isomorphism $V_\ell(\Jac(X/H))\cong V_\ell(\Jac(X))^H.$ Taking dimensions proves \textup{(ii)}.
\end{proof}

\subsection{Reduction and automorphisms}\label{ssec:inertia-automorphisms}
In this subsection, let $C/\mathcal K$ have genus $g\ge1$ and potentially good reduction. If $g=1$, assume that $C$ is an elliptic curve. Choose any finite Galois extension $F/\mathcal K$ over which $\Jac(C)$ has good reduction and put $J=I(F/\mathcal K)$. We write $I=\rho(I_{\mathcal K})$ for the image of inertia on $V$. Here $I$ is a quotient of $J$.

The curve $C$ also has good reduction over $F$. This is immediate when $g=1$, and follows from potentially good reduction and \cite[Thm.~2.4 and Lem.~1.12]{DeligneMumford} when $g\ge2$. Let $\mathcal C/\mathcal O_F$ be its smooth proper model and $\overline C=\mathcal C\times_{\mathcal O_F}\overline k$ its geometric special fibre. This is a smooth projective curve of genus $g$ over $\overline k$.

The finite group $J$ acts on $\overline C$ via $\overline k$-automorphisms. To describe this action, let $s:C_F\xrightarrow{\sim}\mathcal C_F$ be the change of variables to the good-reduction model. Write $F^{\mathrm{nr}}$ for the maximal unramified extension of $F$ in $\overline{\mathcal K}$, and put $R=\mathcal O_{F^{\mathrm{nr}}}$. For $\sigma\in I_{\mathcal K}$, the action is given by the following ``lift--act--reduce" construction \cite[Thm.~1.5]{DDM}:
\begin{equation}\label{eq:lift-act-reduce}
\overline C(\overline k)
\xrightarrow{\mathrm{lift}}\mathcal C(R)
\xrightarrow{s^{-1}}C(F^{\mathrm{nr}})
\xrightarrow{\sigma}C(F^{\mathrm{nr}})
\xrightarrow{s}\mathcal C(R)
\xrightarrow{\mathrm{reduce}}\overline C(\overline k).
\end{equation}
The construction is independent of the chosen lift. Since $I_{\mathcal K}$ fixes $\overline k$ and $I_F$ fixes $F^{\mathrm{nr}}$, it induces an action of $J$ on $\overline C$ by $\overline k$-automorphisms.

\begin{proposition}\label{prop:inertia-automorphisms}
There is an isomorphism of $J$-representations
\begin{equation*}
V\cong V_\ell(\Jac(\overline C))\otimes_{\Q_\ell}\overline{\Q}_\ell,
\end{equation*}
where the action on the right is induced by \eqref{eq:lift-act-reduce}. Moreover, the $J$-action on $\overline C$ factors through a faithful action of $I$, giving an embedding $I\hookrightarrow\Aut_{\overline k}(\overline C)$.
\end{proposition}

\begin{proof}
By \cite[Lem.~4.2(2)]{DDM}, reduction induces a $J$-equivariant isomorphism $T_\ell(\Jac(C))^{I_F}\cong T_\ell(\Jac(\overline C)),$ where the action on the left factors through $I_{\mathcal K}/I_F\cong J$. Since $\Jac(C)$ has good reduction over $F$, the subgroup $I_F$ acts trivially on its Tate module. Tensoring with $\overline{\Q}_\ell$ gives the asserted isomorphism of $J$-representations.

Let $\sigma\in J$ act trivially on $V$, and let $\alpha$ be the automorphism of $\overline C$ induced by \eqref{eq:lift-act-reduce}. By the above isomorphism, $\alpha$ acts trivially on $V_\ell(\Jac(\overline C))$. Suppose first that $g\ge2$. By Lemma \ref{lem:quotient-invariants}(ii), we have $g(\overline C/\langle\alpha\rangle)=g$. The Riemann--Hurwitz formula gives $2g-2\ge|\langle\alpha\rangle|(2g-2)$, so $\alpha=1$ since $g\ge2$.

If $g=1$, then $\alpha$ induces the identity on $\Jac(\overline C)$, so it is either a non-trivial translation or the identity. Since the action in \eqref{eq:lift-act-reduce} fixes the origin of $\overline C$, we deduce that $\alpha=1$.

Conversely, any element of $J$ acting trivially on $\overline C$ acts trivially on $V$. Thus the actions of $J$ on $\overline C$ and on $V$ have the same kernel. Since $I$ is the image of $J$ on $V$, the action on $\overline C$ factors through a faithful action of $I$, as required.
\end{proof}

\begin{remark}\label{rem:J=I}
If $F/\mathcal K$ is chosen to have minimal degree among finite Galois extensions over which $\Jac(C)$ has good reduction, then the natural map $J\to I$ is an isomorphism, so we may identify $J\cong I$. We have stated Proposition~\ref{prop:inertia-automorphisms} for an arbitrary such extension $F$ because, in our later constructions (see, for example, Lemma~\ref{lem:wild-single-reduction}), we choose explicit good-reduction extensions without first establishing minimality. To identify $J$ with $I$, it then suffices to show that the induced $J$-action on $\overline C$ is faithful, as we verify in the relevant constructions below.
\end{remark}

The following example applies Proposition~\ref{prop:inertia-automorphisms} to compute the inertia action on the special fibre from a change of coordinates. In later applications, we give the substitutions and the inertia actions without the intermediate calculations. For further examples, see \cite[\S6]{DD} and \cite[Example~1.9]{DDM}.

\begin{example}\label{ex:inertia_action_tame}
Let $K/\Q_p$ be a finite extension with $p>2$, let $r\ne p$ be an odd prime, and consider the curve $H_r/K$ defined by $y^2=x^r-\pi$. Choose $t\in\overline K$ with $t^{2r}=\pi$, and put $F=K(t,\mu_{2r})$. The substitution $x=t^2u$ and $y=t^rw$ gives a model $\mathcal H_r/\mathcal O_F$ with good reduction and geometric special fibre $\overline H_r:w^2=u^r-1$. Write $s(x,y)=(x/t^2,y/t^r)$ for the change of coordinates.

The group $J=I(F/K)$ is cyclic of order $2r$. Choose $\tau\in I_K$ whose image generates $J$, with $\tau(t)=\zeta t$ for a primitive $2r$th root of unity $\zeta$. Put $R=\mathcal O_{F^{\mathrm{nr}}}$. For an affine point $(u,w)\in\overline H_r(\overline k)$, choose a lift $(U,W)\in\mathcal H_r(R)$. The construction in \eqref{eq:lift-act-reduce} is then
\begin{equation*}
\begin{aligned}
(u,w)
\xrightarrow{\mathrm{lift}}(U,W)
\xrightarrow{s^{-1}}(t^2U,t^rW)
\xrightarrow{\tau}(\zeta^2t^2\tau(U),\zeta^rt^r\tau(W))
\xrightarrow{s}(\zeta^2\tau(U),\zeta^r\tau(W))
\xrightarrow{\mathrm{reduce}}(\overline\zeta^2u,\overline\zeta^rw).
\end{aligned}
\end{equation*}
The last step uses that inertia acts trivially on $\overline k$. Thus $\tau$ induces the automorphism $(u,w)\mapsto(\overline\zeta^2u,\overline\zeta^rw)$ of $\overline H_r$. Since $\overline\zeta^2$ has order $r$ and $\overline\zeta^r=-1$, this automorphism has order $2r$. The $J$-action on $\overline H_r$ is therefore faithful, so Proposition~\ref{prop:inertia-automorphisms} identifies the image of inertia on $V_\ell(\Jac(H_r))$ with $J\cong C_{2r}$.
\end{example}

\subsection{Cluster pictures and Tate modules}\label{ssec:hyperelliptic-clusters}
Throughout this subsection, let $K/\Q_p$ be a finite extension with $p>2$, and let $\ell\ne p$ be prime. We consider hyperelliptic curves of genus at least two whose Jacobians have potentially good reduction. Under this assumption, the description of the Tate module in terms of cluster pictures simplifies.

We let $C/K$ be a hyperelliptic curve of genus $g\ge2$ and choose an equation $y^2=f(x)$ where $f\in K[x]$ is square-free of degree $2g+1$ or $2g+2$. Write $\mathcal R\subset\overline K$ for its set of roots, and extend the valuation $v$ of $K$ to $\overline K$.

A \emph{cluster} is a non-empty subset $\mathfrak s\subseteq\mathcal R$ of the form $\mathfrak s=\{r\in\mathcal R:v(r-z)\ge d\}$ for some $z\in\overline K$ and $d\in\Q$ \cite[Def.~1.1]{DDMM}. We use the following terminology \cite[Defs.~1.1, 1.4 and 1.5]{DDMM}:
\begin{itemize}
\item A cluster is \emph{proper} if it contains more than one root. Its \emph{depth} is $d_{\mathfrak s}=\min_{r\ne r'\in\mathfrak s}v(r-r')$.

\item The \emph{children} of a cluster are the maximal clusters strictly contained in it.

\item A cluster is \emph{odd} or \emph{even} according to its cardinality.

\item A cluster $\mathfrak s$ is \emph{principal} if $|\mathfrak s|\ge3$, with two exceptions: $\mathfrak s=\mathcal R$ is even with exactly two children, or $\mathfrak s$ has a child of cardinality $2g$.
\end{itemize}

Following the conventions of \cite{DDMM}, we draw cluster pictures using dots to represent roots and ovals to represent proper clusters. The children are the dots and ovals immediately inside a given oval. The outermost cluster is labelled by its depth. All other proper clusters are labelled by their relative depths $d_{\mathfrak s}-d_{\mathfrak p}$, where $\mathfrak p$ is the parent.

\begin{example}\label{ex:cluster-decomposition}
Consider the genus-two curve
\begin{equation*}
C/\Q_p:y^2=(x^3-p^3)(x-1)(x+2),
\end{equation*}
where $p>3$. If $\zeta_3$ is a primitive cube root of unity, then $\mathcal R=\{p,p\zeta_3,p\zeta_3^2,1,-2\}$. The roots in $\mathfrak s=\{p,p\zeta_3,p\zeta_3^2\}$ have pairwise differences of valuation $1$, whereas every other pair of distinct roots has difference of valuation $0$. Thus the only proper clusters are $\mathfrak s$ and $\mathcal R$, with depths $1$ and $0$. This gives the cluster picture
\begin{center}
\scalebox{1.8}{\clusterpicture
\Root{1}{first}{r1}
\Root{}{r1}{r2}
\Root{}{r2}{r3}
\ClusterLDName s[][1][\mathfrak s] = (r1)(r2)(r3);
\Root{3}{s}{r4}
\Root{}{r4}{r5}
\ClusterLDName R[][0][\mathcal R] = (s)(sn)(r4)(r5);
\endclusterpicture}
\end{center}
\medskip
Both $\mathfrak s$ and $\mathcal R$ are principal. The children of $\mathcal R$ are $\mathfrak s$, $\{1\}$ and $\{-2\}$, while the children of $\mathfrak s$ are its three singleton roots. We will return to this example after describing the associated Tate-module decomposition.
\end{example}
We will repeatedly use the following reduction criteria \cite[Thm.~1.9(6), (3)]{DDMM}.
\begin{itemize}
\item $\Jac(C)$ has potentially good reduction if and only if every cluster strictly contained in $\mathcal R$ is odd.
\item $C$ has potentially good reduction if and only if every proper cluster contains at least $2g+1$ roots.
\end{itemize}

To describe the Tate module in terms of clusters, we use the following notation.

\begin{notation}\label{not:clusters}
Write $\mathscr P$ for the principal clusters and $\mathcal P=G_K\backslash\mathscr P$. For $\mathfrak s\in\mathscr P$, put
\begin{equation*}
G_{\mathfrak s}=\operatorname{Stab}_{G_K}(\mathfrak s),
\qquad I_{\mathfrak s}=G_{\mathfrak s}\cap I_K,
\qquad P_{\mathfrak s}=G_{\mathfrak s}\cap P_K.
\end{equation*}
Let $\widetilde{\mathfrak s}$ denote the set of odd children of $\mathfrak s$, and put $U_{\mathfrak s}=\textup{ker}(\overline{\Q}_\ell[\tilde{\mathfrak s}]\to\overline{\Q}_\ell),$ where the map sums the coefficients. This is the augmentation representation, which we will also denote by $\overline{\Q}_\ell[\widetilde{\mathfrak s}]\ominus\mathds{1}$. For each principal cluster $\mathfrak s$, let $\widetilde\Gamma_{\mathfrak s}$ be the normalisation of the component $\Gamma_{\mathfrak s}$ defined in \cite[Def.~8.5]{DDMM}. Put $W_{\mathfrak s}=V_\ell(\Jac(\widetilde\Gamma_{\mathfrak s}))\otimes_{\Q_\ell}\overline{\Q}_\ell,$ viewed as a $G_{\mathfrak s}$-representation via the action described in \cite[Thm.~8.8(i)]{DDMM}.
\end{notation}

The following lemma specialises \cite[Thm.~10.1]{DDMM} to the case of potentially good reduction.

\begin{lemma}\label{lem:cluster-tate}
Suppose that $C$ is a hyperelliptic curve of genus $g\ge2$ and $\Jac(C)$ has potentially good reduction. Then the following hold.
\begin{itemize}
\item[\textup{(i)}] There is an isomorphism of $G_K$-representations
\begin{equation*}
V\cong\bigoplus_{[\mathfrak s]\in\mathcal P}\Ind_{G_{\mathfrak s}}^{G_K}W_{\mathfrak s},
\end{equation*}
where each summand is independent of the chosen orbit representative.

\item[\textup{(ii)}] For every $\mathfrak s\in\mathscr P$, there exist a character $\theta_{\mathfrak s}:I_{\mathfrak s}\to\overline{\Q}_\ell^\times$ of finite order prime to $p$ and a representation $\epsilon_{\mathfrak s}$ of $I_{\mathfrak s}$ such that
\begin{equation*}
W_{\mathfrak s}|_{I_{\mathfrak s}}\oplus\epsilon_{\mathfrak s}
\cong\theta_{\mathfrak s}\otimes U_{\mathfrak s}
\end{equation*}
as $I_{\mathfrak s}$-representations. Here $\epsilon_{\mathfrak s}$ is zero if $\mathfrak s$ is odd, and is a one-dimensional character of order at most two if $\mathfrak s$ is even.

\item[\textup{(iii)}] The character of $W_{\mathfrak s}|_{I_{\mathfrak s}}$ is rational for every $\mathfrak s\in\mathscr P$.

\item[\textup{(iv)}] Every child of a proper cluster is odd, and for every $\mathfrak s\in\mathscr P$,
\begin{equation*}
\dim W_{\mathfrak s}=
\begin{cases}
|\widetilde{\mathfrak s}|-1,&\mathfrak s\text{ odd},\\
|\widetilde{\mathfrak s}|-2,&\mathfrak s\text{ even}.
\end{cases}
\end{equation*}
\end{itemize}
\end{lemma}

\begin{proof}
By \cite[Thm.~1.9(6)]{DDMM}, every cluster strictly contained in $\mathcal R$ is odd. Hence every child of a proper cluster is odd, and no principal cluster is \"ubereven, that is, has only even children.

Since $\Jac(C)$ has potentially good reduction, the toric contribution in $H^1$ vanishes. Parts \textup{(i)} and \textup{(ii)} therefore follow from \cite[Thm.~10.1]{DDMM}, after dualising. Here $U_{\mathfrak s}$ and $\epsilon_{\mathfrak s}$ are self-dual, and $\theta_{\mathfrak s}$ is the inverse of the character $\tilde{\gamma}_{\mathfrak s}$ appearing there.

Part \textup{(iii)} follows from Lemma~\ref{lem:quotient-invariants}(i), since $I_{\mathfrak s}$ acts through a finite group of automorphisms on $\widetilde\Gamma_{\mathfrak s}$. Finally, the dimension formula in \textup{(iv)} follows from \textup{(ii)}, since $\dim U_{\mathfrak s}=|\widetilde{\mathfrak s}|-1$ and $\dim\epsilon_{\mathfrak s}$ is zero or one according as $\mathfrak s$ is odd or even.
\end{proof}

\begin{example}
We return to Example~\ref{ex:cluster-decomposition}. Every cluster strictly contained in $\mathcal R$ is odd, so $\Jac(C)$ has potentially good reduction. Both principal clusters are fixed by $G_{\Q_p}$ and have three odd children. Lemma~\ref{lem:cluster-tate} therefore gives
\begin{equation*}
V\cong W_{\mathcal R}\oplus W_{\mathfrak s},
\qquad
\dim W_{\mathcal R}=\dim W_{\mathfrak s}=2.
\end{equation*}
Thus the Tate module is not absolutely irreducible.
\end{example}

Lemma~\ref{lem:cluster-tate} expresses the contribution of each principal cluster in terms of the permutation representation on its children. We therefore describe this action explicitly. In suitable coordinates, inertia acts by affine transformations and wild inertia acts by translations. To introduce these coordinates, fix a principal cluster $\mathfrak s$, choose $z\in\mathfrak s$ and $t\in\overline K$ with $v(t)=d_{\mathfrak s}$, and define
\begin{equation}\label{eq:injection_roots}
\iota_{\mathfrak s}:\mathfrak s\to\overline k,
\qquad
r\mapsto\frac{r-z}{t}\mod{\mathfrak m_{\overline K}}.
\end{equation}

\begin{lemma}\label{lem:cluster-inertia}
For every principal cluster $\mathfrak s$, the following hold.
\begin{itemize}
\item[\textup{(i)}] Two roots of $\mathfrak s$ have the same image under $\iota_{\mathfrak s}$ if and only if they belong to the same child.

\item[\textup{(ii)}] For $h\in I_{\mathfrak s}$ and $r\in\mathfrak s$, we have
\begin{equation*}
\iota_{\mathfrak s}(h(r))=\alpha_{\mathfrak s}(h)\iota_{\mathfrak s}(r)+\beta_{\mathfrak s}(h),
\end{equation*}
where $\alpha_{\mathfrak s}(h)=h(t)/t\mod{\mathfrak m_{\overline K}}$ and $\beta_{\mathfrak s}(h)=(h(z)-z)/t\mod{\mathfrak m_{\overline K}}$.

\item[\textup{(iii)}] For $h\in P_{\mathfrak s}$ and $r\in\mathfrak s$, we have
\begin{equation*}
\iota_{\mathfrak s}(h(r))=\iota_{\mathfrak s}(r)+\beta_{\mathfrak s}(h).
\end{equation*}

\item[\textup{(iv)}] All $P_{\mathfrak s}$-orbits on the children of $\mathfrak s$ have the same cardinality $p^{a_{\mathfrak s}}$ for some integer $a_{\mathfrak s}\ge0$. In particular, $p^{a_{\mathfrak s}}\mid|\mathfrak s|$.
\end{itemize}
\end{lemma}

\begin{proof}
The expression defining $\iota_{\mathfrak s}$ is integral because $v(r-z)\ge d_{\mathfrak s}$. Two roots $r,r'\in\mathfrak s$ have the same reduction precisely when $v(r-r')>d_{\mathfrak s}$. This is equivalent to $r$ and $r'$ belonging to the same child. This proves (i).

Since $h\in I_{\mathfrak s}$ fixes $\overline k$ and stabilises $\mathfrak s$, reducing the identity
\begin{equation*}
\frac{h(r)-z}{t}
=\frac{h(t)}{t}h\left(\frac{r-z}{t}\right)+\frac{h(z)-z}{t}
\end{equation*}
gives the asserted formula in part (ii). The map $\alpha_{\mathfrak s}:I_{\mathfrak s}\to\overline k^\times$ is a character with finite image and of order prime to $p$. It is therefore trivial on $P_{\mathfrak s}$; cf.~\cite[Rem.~8.3]{DDMM}. This proves (iii).

By \textup{(i)} and \textup{(iii)}, an element $h\in P_{\mathfrak s}$ fixes a child if and only if $\beta_{\mathfrak s}(h)=0$, in which case it fixes every child. Thus all children have the same stabiliser in $P_{\mathfrak s}$, so their orbits have equal cardinality. This cardinality is a power of $p$, since $P_{\mathfrak s}$ is pro-$p$. Children in the same orbit contain equally many roots, so each orbit gives a multiple of $p^{a_{\mathfrak s}}$ roots. As the children partition $\mathfrak s$, we deduce $p^{a_{\mathfrak s}}\mid|\mathfrak s|$.
\end{proof}

\section{Hyperelliptic curves with absolutely irreducible Tate modules}\label{sec:hyperelliptic-irreducibility}

This section and the next prove Theorem~\ref{thm:hyperelliptic-classification}. Throughout, let $K/\Q_p$ be a finite extension with $p>2$, let $q$ be the cardinality of its residue field, and let $\ell\ne p$ be prime.

A curve $C/K$ of genus at least two is \emph{geometrically hyperelliptic} if $C_{\overline K}$ admits a morphism of degree two to $\mathbb P^1_{\overline K}$. We call $C$ \emph{hyperelliptic over $K$}, or simply \emph{hyperelliptic}, if it admits a morphism of degree two to $\mathbb P^1_K$ over $K$. For a geometrically hyperelliptic curve, this is equivalent to its quotient by the hyperelliptic involution having a $K$-point.

In \S\ref{ssec:no-K-equation}, we show that every geometrically hyperelliptic curve $C/K$ of genus at least two for which $V_\ell(\Jac(C))$ is absolutely irreducible is hyperelliptic over $K$. The classifications below therefore also apply to geometrically hyperelliptic curves.

We begin with genus one, where the necessary conditions follow directly from the dimension restrictions.

\begin{proposition}\label{prop:elliptic-classification}
There exists a genus-one curve $C/K$ such that $V_\ell(\Jac(C))$ is absolutely irreducible and
\begin{itemize}
\item[\textup{(i)}] tamely ramified if and only if $q\equiv2\pmod3$ or $q\equiv3\pmod4$;
\item[\textup{(ii)}] wildly ramified if and only if $p=3$.
\end{itemize}
\end{proposition}

\begin{proof}
It suffices to consider elliptic curves, since the Jacobians of genus-one curves are elliptic curves. Necessity follows from Lemma~\ref{lem:tame}, which applies by Lemma~\ref{lem:arithmetic-consequences}\textup{(iii)}, while sufficiency follows from the constructions below.

In the tame case, Lemma~\ref{lem:tame}\textup{(i)} requires $\varphi(m)=2$, so $m\in\{3,4,6\}$. The additional conditions $(m,p)=1$ and $\ord_m(q)=2$ hold for $m=3$ precisely when $q\equiv2\pmod3$, for $m=4$ precisely when $q\equiv3\pmod4$, and, since $q$ is odd, for $m=6$ precisely when $q\equiv2\pmod3$. Proposition~\ref{prop:tame-superelliptic-realisation} gives elliptic curves with tamely ramified, absolutely irreducible Tate modules by taking $(\nu,a)=(3,1)$ or $(2,2)$, respectively.

In the wild case, Lemma~\ref{lem:tame}\textup{(ii)} requires $(p-1)\mid2$, so $p=3$. Conversely, Proposition~\ref{prop:wild-single-cluster}, with $p=3$, $a=1$ and $b=\eta=1$, shows that the elliptic curve $E:y^2=x^3-\pi x-\pi$ has a wildly ramified, absolutely irreducible Tate module.
\end{proof}
For genus at least two, the tame and wild classifications follow by combining the restrictions proved later in this section and the constructions in the next.

\begin{theorem}\label{thm:tamehyperelliptic}
Let $g\ge2$. There exists a hyperelliptic curve $C/K$ of genus $g$ whose Tate module $V_\ell(\Jac(C))$ is absolutely irreducible and tamely ramified if and only if $2g+1$ is a prime distinct from $p$ with $\ord_{2g+1}(q)=2g$.
\end{theorem}

\begin{proof}
Necessity follows from Corollary~\ref{cor:tame-hyperelliptic-genus}, and sufficiency from Corollary~\ref{cor:tame-hyperelliptic-realisation}.
\end{proof}

\begin{theorem}\label{thm:wildhyperelliptic}
Let $g\ge2$. There exists a hyperelliptic curve $C/K$ of genus $g$ whose Tate module $V_\ell(\Jac(C))$ is absolutely irreducible and wildly ramified if and only if $g={p^a-1}/{2}$ or $g=p^a-1$ for some integer $a\ge1$.
\end{theorem}

\begin{proof}
Necessity follows from Corollary~\ref{cor:wild-hyperelliptic-genus}, and sufficiency from Propositions~\ref{prop:wild-single-cluster} and \ref{prop:wild-two-clusters}.
\end{proof}

Using Lemmas~\ref{lem:cluster-tate} and~\ref{lem:cluster-inertia}, we restrict the cluster picture to the two shapes in Corollary~\ref{cor:cluster-configurations}. We then treat the tame and wild restrictions on the genus separately.

\subsection{Restrictions on the cluster picture}\label{ssec:cluster-restrictions}

For a curve $C/K$, put $V=V_\ell(\Jac(C))\otimes_{\Q_\ell}\overline{\Q}_\ell$.

If $V$ is irreducible, then $\Jac(C)$ has potentially good reduction by Lemma~\ref{lem:arithmetic-consequences}, so Lemma~\ref{lem:cluster-tate} applies. As in Notation~\ref{not:clusters}, we write $\mathscr P$ for the set of principal clusters. Recall that a cluster $\mathfrak s$ is principal if $|\mathfrak s|\ge3$, unless $\mathfrak s=\mathcal R$ is even with exactly two children or $\mathfrak s$ has a child of cardinality $2g$.

\begin{lemma}\label{lem:principal-orbit}
Let $C/K$ be a hyperelliptic curve of genus $g\ge2$. Suppose that $V_\ell(\Jac(C))$ is absolutely irreducible. Then the following hold.
\begin{itemize}
\item[\textup{(i)}] $G_K$ acts transitively on $\mathscr P$, and $|\mathscr P|\in\{1,2\}$.

\item[\textup{(ii)}] If $\mathscr P$ contains an even cluster, then $\mathscr P=\{\mathcal R\}$.

\item[\textup{(iii)}] Every child of a principal cluster is a singleton.
\end{itemize}
\end{lemma}

\begin{proof}
Since $\Jac(C)$ has potentially good reduction, every cluster strictly contained in $\mathcal R$ is odd. Hence, every non-singleton cluster strictly contained in $\mathcal R$ is principal.

Every odd principal cluster has at least three children. An even principal cluster must be $\mathcal R$ and has at least four children. Indeed, all children of $\mathcal R$ are odd and $\mathcal R$ cannot have exactly two children as it is assumed to be principal. Lemma~\ref{lem:cluster-tate}\textup{(iv)} therefore gives $\dim W_{\mathfrak s}\ge2$ for every $\mathfrak s\in\mathscr P$. Consequently, every summand in Lemma~\ref{lem:cluster-tate}\textup{(i)} is non-zero. Since $V$ is irreducible, this decomposition has exactly one summand. Hence $G_K$ acts transitively on $\mathscr P$. This completes the first half of part (i).

Since $G_K$ acts transitively on $\mathscr P$, all principal clusters have the same cardinality. If $\mathcal R$ is principal, it is fixed by $G_K$, so $\mathscr P=\{\mathcal R\}$. Otherwise, $\mathcal R$ has exactly two children, since it cannot have a child of cardinality $2g$. Each non-singleton child is principal and cannot properly contain another principal cluster, because the two would have different cardinalities. Thus $\mathscr P$ consists precisely of the non-singleton children of $\mathcal R$, giving $|\mathscr P|\in\{1,2\}$. This completes the proof of (i).

For \textup{(ii)}, an even principal cluster must be $\mathcal R$. Since $\mathcal R$ is fixed by $G_K$, transitivity gives $\mathscr P=\{\mathcal R\}$. For \textup{(iii)}, a non-singleton child of a principal cluster would itself be principal, contradicting the fact that all principal clusters have the same cardinality.
\end{proof}

\begin{corollary}\label{cor:cluster-configurations}
Let $C/K$ be a hyperelliptic curve of genus $g\ge2$. Suppose that $V_\ell(\Jac(C))$ is absolutely irreducible. Then $C$ admits an affine equation $y^2=f(x)$, where $f\in K[x]$ is square-free of degree $2g+1$ or $2g+2$, and its cluster picture has one of the following two shapes:
\begin{center}
\begin{tabular}{cc}
\scalebox{1.8}{\clusterpicture
  \Root{1}{first}{r1}
  \Root{}{r1}{r2}
  \RootT{}{r2}{dots}{$\cdots$}
  \Root{}{dots}{r3}
  \ClusterLDName R[][][\mathcal R] = (r1)(r2)(dots)(r3);
\endclusterpicture}
&\quad\quad\quad
\scalebox{1.8}{\clusterpicture
  \Root{1}{first}{r1}
  \Root{}{r1}{r2}
  \RootT{}{r2}{dots1}{$\cdots$}
  \Root{}{dots1}{r3}
  \ClusterLDName s1[][][\mathfrak s_1] = (r1)(r2)(dots1)(r3);
  \Root{3}{s1}{r4}
  \Root{}{r4}{r5}
  \RootT{}{r5}{dots2}{$\cdots$}
  \Root{}{dots2}{r6}
  \ClusterLDName s2[][][\mathfrak s_2] = (r4)(r5)(dots2)(r6);
  \ClusterLDName R[][][\mathcal R] = (s1)(s1n)(s2)(s2n);
\endclusterpicture}
\\[6pt]
\textup{(i)} & \textup{(ii)}
\\[3pt]
$\mathscr P=\{\mathcal R\}$
&$\mathscr P=\{\mathfrak s_1,\mathfrak s_2\}$
\end{tabular}
\end{center}
\begin{itemize}
\item[\textup{(i)}] $\mathcal R$ is the only proper cluster and may be odd or even.

\item[\textup{(ii)}] $\mathcal R=\mathfrak s_1\sqcup\mathfrak s_2$, where $\mathfrak s_1$ and $\mathfrak s_2$ are $G_K$-conjugate clusters, each containing $g+1$ roots. Their cardinality is odd, so $g$ is even.
\end{itemize}
\end{corollary}

\begin{proof}
Choose an equation $C:y^2=f(x)$ with $\deg f\in\{2g+1,2g+2\}$. As in the proof of \hbox{Lemma \ref{lem:principal-orbit}}, every cluster strictly contained in $\mathcal R$ is odd. If $\mathcal R$ is principal, all its children are singletons by Lemma~\ref{lem:principal-orbit}\textup{(iii)}. This gives picture \textup{(i)}.

Suppose instead that $\mathcal R$ is not principal. By definition, either it is even with exactly two children, or it has a child of cardinality $2g$. The latter is impossible because every child is odd. Thus $\mathcal R$ has exactly two children and $|\mathcal R|=2g+2$. Each non-singleton child is odd of cardinality at least three and hence principal. The two children cannot both be singletons, since $|\mathcal R|\ge6$.

If exactly one child is a singleton, the other is a principal cluster $\mathfrak s$ containing $2g+1$ roots. Its children are all singletons, giving
\begin{center}
\scalebox{1.6}{\clusterpicture
  \Root{1}{first}{r1}
  \Root{}{r1}{r2}
  \RootT{}{r2}{dots}{$\cdots$}
  \Root{}{dots}{r3}
  \ClusterLDName s[][][\mathfrak s] = (r1)(r2)(dots)(r3);
  \Root{3}{s}{r4}
  \ClusterLDName R[][][{}] = (s)(sn)(r4);
\endclusterpicture}
\end{center}
Here $\mathcal R=\mathfrak s\sqcup\{r_0\}$ and $|\mathfrak s|=2g+1$. Since the two children have different cardinalities, $G_K$ fixes each of them, so $r_0\in K$. Translating $r_0$ to $0$ gives an equation $y^2=xh(x)$ with $\deg h=2g+1$ and $h(0)\ne0$. The substitution $u=1/x$, $w=y/x^{g+1}$ gives a $K$-isomorphic curve with equation
\begin{equation*}
w^2=u^{2g+1}h(1/u).
\end{equation*}
Its defining polynomial is square-free of degree $2g+1$, with roots $1/r$ for $r\in\mathfrak s$.

Since $\mathfrak s$ and $\{0\}$ are the two children of $\mathcal R$, we have $v(r)=d_{\mathcal R}$ for every $r\in\mathfrak s$. Since the children of $\mathfrak s$ are singletons, we also have $v(r-r')=d_{\mathfrak s}$ for distinct $r,r'\in\mathfrak s$. Consequently,
\begin{equation*}
v(1/r-1/r')=d_{\mathfrak s}-2d_{\mathcal R}
\end{equation*}
for every such pair. Thus the transformed roots form the only proper cluster, giving picture \textup{(i)}.

If both children are non-singletons, they are principal and have only singleton children, giving picture \textup{(ii)}. By Lemma~\ref{lem:principal-orbit}\textup{(i)}, they are Galois-conjugate and therefore have equal cardinality $g+1$. Since both are odd, $g$ is even.
\end{proof}

\subsection{The tame case}
In this subsection, $C/K$ is a hyperelliptic curve of genus $g\ge2$ with an absolutely irreducible and tamely ramified Tate module $V$. We write $I$ for its finite inertia image and $m=|I|$. By Lemma~\ref{lem:tame}\textup{(i)}, the group $I$ is cyclic of order prime to $p$, and
\begin{equation}\label{eq:tame-hyperelliptic-dimension}
2g=\varphi(m)=\ord_m(q).
\end{equation}
Since $q$ generates $(\Z/m\Z)^\times$ and $g\ge2$, it follows that
\begin{equation}\label{eq:tame-hyperelliptic-order}
m=\nu^a\quad\text{or}\quad m=2\nu^a
\end{equation}
for an odd prime $\nu\ne p$ and an integer $a\ge1$.

We show that $a=1$, so that $2g+1=\nu$ is prime and $\ord_\nu(q)=\nu-1$. This gives the necessary condition in Theorem~\ref{thm:tamehyperelliptic}.
\begin{lemma}\label{lem:tame-principal}
The set $\mathscr P$ consists of a single cluster. Consequently, $C$ has potentially good reduction.
\end{lemma}

\begin{proof}
By Lemma~\ref{lem:principal-orbit}\textup{(i)}, it suffices to exclude $\mathscr P=\{\mathfrak s_1,\mathfrak s_2\}$. Suppose this occurs. Lemma~\ref{lem:cluster-tate}\textup{(i)} gives $V=W_{\mathfrak s_1}\oplus W_{\mathfrak s_2},$ where both summands are non-zero and $G_K$ permutes them according to its action on the clusters.

The action of $I_K$ on the two principal clusters factors through $I$. Indeed, an element acting trivially on $V$ cannot exchange the two non-zero summands.

Choose a generator $\tau$ of the cyclic group $I$. Its eigenvalues on $V$ are the primitive $m$th roots of unity, each occurring once, by \cite[proof of Prop.~3]{ACampo}.

Suppose first that $I$ stabilises both clusters. Each $W_{\mathfrak s_i}$ then has rational character by Lemma~\ref{lem:cluster-tate}\textup{(iii)}, so the characteristic polynomial $P_i(T)$ of $\tau$ on $W_{\mathfrak s_i}$ lies in $\Q[T]$. The decomposition above gives
\begin{equation*}
\Phi_m(T)=P_1(T)P_2(T),
\end{equation*}
where $\Phi_m$ is the $m$th cyclotomic polynomial. Both factors have positive degree, contradicting the irreducibility of $\Phi_m$ over $\Q$.

Suppose instead that $I$ exchanges the two clusters. Then $\tau$ exchanges the two summands, so $m$ is even. By \eqref{eq:tame-hyperelliptic-order}, we have $m=2\nu^a$ with $\nu$ odd. Since every eigenvalue of $\tau$ is a primitive $2\nu^a$th root of unity, we have $\tau^{\nu^a}=-\operatorname{id}_V$, which preserves both summands. However, $\nu^a$ is odd, so $\tau^{\nu^a}$ exchanges them. This is again a contradiction.

Thus $\mathscr P$ consists of a single cluster. By Corollary~\ref{cor:cluster-configurations}, we may choose an equation whose only proper cluster is $\mathcal R$. The curve $C$ therefore has potentially good reduction.
\end{proof}

We may now study inertia through its faithful action on the smooth special fibre. Its action on the branch locus of the hyperelliptic cover forces $a=1$ in \eqref{eq:tame-hyperelliptic-order}.

\begin{lemma}\label{lem:tame-hyperelliptic-inertia}
The finite group $I$ is cyclic of order $\nu$ or $2\nu$ for an odd prime $\nu\ne p$.
\end{lemma}

\begin{proof}
By Lemma~\ref{lem:tame-principal}, the curve $C$ acquires good reduction, and its geometric special fibre $\overline C$ is hyperelliptic by \cite[Thm.~8.8]{DDMM}. By Proposition~\ref{prop:inertia-automorphisms}, the group $I$ acts faithfully on $\overline C$. Its order is $\nu^a$ or $2\nu^a$, where $\nu\ne p$ is an odd prime, by \eqref{eq:tame-hyperelliptic-order}. Let $\iota$ be the hyperelliptic involution of $\overline C$, and let $B$ be the branch locus of the hyperelliptic cover $\overline C\to\overline C/\langle\iota\rangle=\mathbb P^1_{\overline k}$. Every automorphism of $\overline C$ descends to an automorphism of $\mathbb P^1_{\overline k}$ preserving $B$. Since $|B|=2g+2\ge6$, an automorphism of $\mathbb P^1_{\overline k}$ fixing every point of $B$ is the identity. Hence the kernel of the action of $\Aut_{\overline k}(\overline C)$ on $B$ is $\langle\iota\rangle$.

Let $H\subseteq I$ be the subgroup of order $\nu^a$. Since $H$ has odd order, it acts faithfully on $B$. A generator of $H$ therefore induces a permutation of order $\nu^a$. Such a permutation has a cycle of length $\nu^a$ which gives $\nu^a\le|B|=2g+2=\nu^{a-1}(\nu-1)+2.$
Thus $\nu^{a-1}\le2$, and since $\nu$ is odd, we deduce that $a=1$.
\end{proof}

\begin{corollary}\label{cor:tame-hyperelliptic-genus}
Let $C/K$ be a hyperelliptic curve of genus $g\ge2$. Suppose that $V_\ell(\Jac(C))$ is absolutely irreducible and tamely ramified. Then $2g+1$ is a prime distinct from $p$, and $\ord_{2g+1}(q)=2g$. Moreover, $C$ has potentially good reduction.
\end{corollary}

\begin{proof}
By Lemma~\ref{lem:tame-hyperelliptic-inertia}, we have $m=\nu$ or $2\nu$ for an odd prime $\nu\ne p$. In either case, $\varphi(m)=\nu-1$ and $\ord_m(q)=\ord_\nu(q)$, since $q$ is odd. Equation~\eqref{eq:tame-hyperelliptic-dimension} therefore gives $2g=\nu-1=\ord_\nu(q)$. The final assertion follows from Lemma~\ref{lem:tame-principal}.
\end{proof}

\subsection{The wild case}\label{ssec:wild-case}

In this subsection, let $C/K$ be a hyperelliptic curve of genus $g\ge2$ whose Tate module is absolutely irreducible and wildly ramified. Using the cluster pictures in Corollary~\ref{cor:cluster-configurations}, we determine the number of roots in each principal cluster and hence the possible genera.

\begin{lemma}\label{lem:wild-children}
Suppose that $V$ is absolutely irreducible and wildly ramified. Then there exists an integer $a\ge1$ such that the following hold for every $\mathfrak s\in\mathscr P$.
\begin{itemize}
\item[\textup{(i)}] If $\mathfrak s$ is odd, then $|\mathfrak s|=p^a$.
\item[\textup{(ii)}] If $\mathfrak s$ is even, then $|\mathfrak s|=2p^a$.
\end{itemize}
\end{lemma}

\begin{proof}
By Lemma~\ref{lem:arithmetic-consequences}, $\Jac(C)$ has potentially good reduction, so Lemma~\ref{lem:cluster-tate} applies. Moreover, $V^{P_K}=0$. Indeed, $P_K$ is normal in $G_K$, so $V^{P_K}$ is a $G_K$-stable subspace. Since $V$ is irreducible and wildly ramified, $V^{P_K}=0$.

By Lemma~\ref{lem:principal-orbit}\textup{(i)}, there are at most two principal clusters. Since $P_K$ is pro-$p$ and $p>2$, it stabilises each of them. Using Lemma~\ref{lem:cluster-tate}\textup{(i)}, we deduce that
\begin{equation*}
V|_{P_K}\cong\bigoplus_{\mathfrak s\in\mathscr P}W_{\mathfrak s}|_{P_K}.
\end{equation*}
Since $V^{P_K}=0$, we have $W_{\mathfrak s}^{P_K}=0$ for every $\mathfrak s\in\mathscr P$.

Fix $\mathfrak s\in\mathscr P$, and put $e_{\mathfrak s}=1$ if $\mathfrak s$ is odd and $e_{\mathfrak s}=2$ otherwise. By Lemma~\ref{lem:cluster-tate}\textup{(ii)}, $\theta_{\mathfrak s}$ has order prime to $p$, while $\epsilon_{\mathfrak s}$ is zero or has order at most two. Since $p>2$, their restrictions satisfy
\begin{equation*}
\theta_{\mathfrak s}|_{P_K}\cong\mathds{1},
\qquad
\epsilon_{\mathfrak s}|_{P_K}\cong\mathds{1}^{\oplus(e_{\mathfrak s}-1)}.
\end{equation*}
Every child of $\mathfrak s$ is a singleton by Lemma~\ref{lem:principal-orbit}\textup{(iii)}, so the same formula gives
\begin{equation}\label{eq:W|PK}
W_{\mathfrak s}|_{P_K}\oplus\mathds{1}^{\oplus e_{\mathfrak s}}\cong\overline{\Q}_\ell[\mathfrak s].
\end{equation}
Taking invariants shows that $P_K$ has exactly $e_{\mathfrak s}$ orbits on $\mathfrak s$. Since $P_{\mathfrak s}=P_K$, Lemma~\ref{lem:cluster-inertia}\textup{(iv)} shows that these orbits all have cardinality $p^{a_{\mathfrak s}}$ for some integer $a_{\mathfrak s}\ge0$. Consequently, $|\mathfrak s|=e_{\mathfrak s}p^{a_{\mathfrak s}}$. As $|\mathfrak s|\ge3$ and $e_{\mathfrak s}\le2$, we have $a_{\mathfrak s}\ge1$.

Finally, the principal clusters are Galois-conjugate by Lemma~\ref{lem:principal-orbit}\textup{(i)} and therefore have the same cardinality. Hence both $e_{\mathfrak s}$ and $a_{\mathfrak s}$ are independent of $\mathfrak s$. Writing $a$ for the common exponent proves both assertions.
\end{proof}

Combining the preceding lemma with the cluster pictures in Corollary~\ref{cor:cluster-configurations} gives the following restriction on the genus.
\begin{corollary}\label{cor:wild-hyperelliptic-genus}
Let $C/K$ be a hyperelliptic curve of genus $g\ge2$. Suppose that $V_\ell(\Jac(C))$ is absolutely irreducible and wildly ramified. Then there exists an integer $a\ge1$ such that $g={p^a-1}/{2}$ or $g=p^a-1.$
\end{corollary}

\begin{proof}
Choose an equation for $C$ as in Corollary~\ref{cor:cluster-configurations}. In picture \textup{(i)}, Lemma~\ref{lem:wild-children} gives $|\mathcal R|=p^a$ or $2p^a$ according as $\mathcal R$ is odd or even. In picture \textup{(ii)}, the two principal clusters each contain $p^a$ roots, so $|\mathcal R|=2p^a$. Since the defining polynomial has degree $2g+1$ or $2g+2$, the odd case gives $2g+1=p^a$, while the even case gives $2g+2=2p^a$. Hence $g=(p^a-1)/2$ or $g=p^a-1$, as required.
\end{proof}

\subsection{Geometrically hyperelliptic curves}\label{ssec:no-K-equation}

We now prove that absolute irreducibility forces a geometrically hyperelliptic curve of genus $g\ge2$ to be hyperelliptic over $K$.

\begin{proposition}\label{prop:equation-over-K}
Let $C/K$ be a geometrically hyperelliptic curve of genus $g\ge2$ such that $V_\ell(\Jac(C))$ is absolutely irreducible. Then $C$ is hyperelliptic over $K$ and admits an equation $y^2=f(x)$ with $f\in K[x]$ square-free of degree $2g+1$ or $2g+2$.
\end{proposition}

\begin{proof}
Let $\iota$ be the hyperelliptic involution, defined over $K$ by uniqueness, and put $Q=C/\langle\iota\rangle$. It suffices to show that $Q(K)\ne\emptyset$. For even $g$, this follows from \cite[Thm.~1(4)]{RoeXarles}, since the $g^1_2$ is unique.

Suppose that $g$ is odd, so $g\ge3$. Since $p\ne2$, the representation $V=V_2(\Jac(C))$ is absolutely irreducible by \cite[Cor.~1.6]{DDWeil}. Write $G,I,P$ for the images of $G_K,I_K,P_K$ on $V$. The group $I$ is finite by Lemma~\ref{lem:arithmetic-consequences}. If $P\ne1$, put $H=P$. Otherwise, $I$ is cyclic with $\varphi(|I|)=2g$ by Lemma~\ref{lem:tame}. Since $g\ge3$ is odd, $2g$ is not a power of two, so neither is $|I|$. We take $H$ to be its non-trivial subgroup of largest odd order. In either case, $H$ is a non-trivial characteristic subgroup of $I$ of odd order. It is therefore normal in $G$, and absolute irreducibility gives $V^H=0$.

Let $B\subset Q(\overline K)$ be the branch locus of $C\to Q$, so that $|B|=2g+2$ by Riemann--Hurwitz. For $b\in B$, let $w_b\in C(\overline K)$ be the Weierstrass point above $b$. Write $\F_2[B]^0\subset\F_2[B]$ for the subspace of elements with coefficient sum zero, which contains $\Sigma_B=\sum_{b\in B}b$. Sending $b+b'$ to the class of $w_b-w_{b'}$ induces a $G_K$-equivariant isomorphism $\F_2[B]^0/\langle\Sigma_B\rangle\cong\Jac(C)[2]$ \cite[Thm.~3.6 and \S3.3.1]{LSTX}. As $|B|\ge6$, the symmetric group on $B$ acts faithfully on $\F_2[B]^0/\langle\Sigma_B\rangle$ \cite[Lem.~2.1]{LSTX}. Hence every element of $G_K$ acting trivially on $V$, and therefore on $\Jac(C)[2]$, fixes $B$ pointwise, so the action on $B$ factors through $G$.

Since $|H|$ is odd, $\F_2[H]$-modules are semisimple, and hence $\F_2[B]\cong\Jac(C)[2]\oplus\mathds{1}^{\oplus2}$ as $H$-modules. For the same reason, taking $H$-invariants is exact on $\Z_2[H]$-modules, so $\Jac(C)[2]^H$ is the reduction of $T_2(\Jac(C))^H$, which has rank $\dim V^H=0$. Therefore
\begin{equation*}
\#(H\backslash B)=\dim_{\F_2}\F_2[B]^H=2+\dim_{\F_2}\Jac(C)[2]^H=2.
\end{equation*}

The two $H$-orbits on $B$ have odd cardinalities whose sum is $2g+2$. These cardinalities are unequal, since otherwise each would equal the even integer $g+1$. Since $H$ is normal in $G$, every element of $G_K$ sends an $H$-orbit to one of the same cardinality, so each orbit is $G_K$-stable. Either orbit therefore defines a $K$-rational divisor of odd degree on $Q$, and hence $Q(K)\ne\emptyset$ by \cite[Thm.~14(2)]{RoeXarles}.
\end{proof}

\section{Realisations by hyperelliptic curves}\label{sec:hyperelliptic-realisations}

We now construct the hyperelliptic curves required to complete the proof of Theorem~\ref{thm:hyperelliptic-classification}. Throughout, let $K/\Q_p$ be a finite extension with $p>2$, and let $\ell\ne p$ be prime. We begin with the tame case.

\subsection{Tame realisations}\label{sec:tame-hyperelliptic-realisations}

The genus-one case is covered by Proposition~\ref{prop:elliptic-classification}. For the remaining cases, we consider the hyperelliptic curves
\begin{equation*}
H_r:y^2=x^r-\pi,
\end{equation*}
where $r\ge5$ is prime and distinct from $p$. Each curve acquires good reduction over a tame extension, allowing us to compute the action of inertia on its Tate module through the smooth special fibre. We write $\mathcal R_r$ for the set of roots of $x^r-\pi$. Since $p\nmid r$, any two distinct roots have difference of valuation $1/r$, giving the cluster picture
\begin{center}
\scalebox{1.8}{\clusterpicture
  \Root{1}{first}{r1}
  \Root{}{r1}{r2}
  \RootT{}{r2}{dots}{$\cdots$}
  \Root{}{dots}{r3}
  \ClusterLDName R[][\frac1r][\mathcal R_r]
    = (r1)(r2)(dots)(r3);
\endclusterpicture}
\end{center}
This is the cluster picture of Corollary~\ref{cor:cluster-configurations}\textup{(i)}. To determine when the Tate module is absolutely irreducible, we use the description of inertia from \S\ref{ssec:inertia-automorphisms}.

\begin{lemma}\label{lem:tame-hyperelliptic-reduction}
Choose $t\in\overline K$ with $t^{2r}=\pi$, and put $F=K(t,\mu_{2r})$. The curve $H_r$ has genus $(r-1)/2$ and good reduction over the tame Galois extension $F/K$.
\end{lemma}

\begin{proof}
Since $x^r-\pi$ is square-free of odd degree $r$, the curve $H_r$ has genus $(r-1)/2$. The field $F$ is the splitting field of $x^{2r}-\pi$, so $F/K$ is Galois. Moreover, $K(\mu_{2r})/K$ is unramified, and $x^{2r}-\pi$ is Eisenstein over $K(\mu_{2r})$. Thus $F/K(\mu_{2r})$ is totally ramified of degree $2r$. Since $p\nmid2r$, the extension $F/K$ is tame. The good reduction of $H_r$ over $F$ was established in Example~\ref{ex:inertia_action_tame}.
\end{proof}

\begin{lemma}\label{lem:tame-hyperelliptic-irreducibility}
The representation $V_\ell(\Jac(H_r))$ is tamely ramified and is absolutely irreducible if and only if $\ord_r(q)=r-1$.
\end{lemma}

\begin{proof}
Put $V=V_\ell(\Jac(H_r))\otimes_{\Q_\ell}\overline{\Q}_\ell$. It is tamely ramified by Lemma~\ref{lem:tame-hyperelliptic-reduction}. By Example~\ref{ex:inertia_action_tame}, we may identify the image $I$ of inertia with $I(F/K)\cong C_{2r}$, acting on the geometric special fibre $\overline H_r:w^2=u^r-1$ as described there. Its subgroup $N$ of order $r$ acts through the automorphisms
\begin{equation*}
\{(u,w)\mapsto(\xi u,w):\xi\in\mu_r(\overline k)\}.
\end{equation*}
The quotient $\overline H_r/N$ is therefore $\mathbb P^1_{\overline k}$ with parameter $w$. Proposition~\ref{prop:inertia-automorphisms} and Lemma~\ref{lem:quotient-invariants} give $V^N=0$ and show that the trace of every element of $N$ on $V$ is rational.

Since $N$ has prime order $r$, its non-trivial characters are Galois conjugate over $\Q$. Rationality of the character therefore forces them to occur with equal multiplicities. As $\dim V=r-1$ and $V^N=0$, each occurs exactly once, which gives
\begin{equation*}
V|_N\cong\bigoplus_{\substack{\chi\in\widehat N\\\chi\ne\mathds{1}}}\chi.
\end{equation*}

Since $I$ is abelian, each $N$-character space is $I$-stable. These one-dimensional $I$-representations are pairwise non-isomorphic, since their restrictions to $N$ are distinct. Lemma~\ref{lem:frobenius-transitivity} therefore shows that $V$ is irreducible if and only if $\ord_r(q)=r-1$.
\end{proof}

\begin{corollary}\label{cor:tame-hyperelliptic-realisation}
Let $g\ge2$ satisfy $\ord_{2g+1}(q)=2g$. Then the hyperelliptic curve
\begin{equation*}
H_{2g+1}:y^2=x^{2g+1}-\pi
\end{equation*}
has genus $g$, and $V_\ell(\Jac(H_{2g+1}))$ is tamely ramified and absolutely irreducible.
\end{corollary}

\begin{proof}
Put $r=2g+1$. Since $r-1=\ord_r(q)\le\varphi(r)\le r-1,$ we have $\varphi(r)=r-1$, so $r$ is prime. Moreover, $(r,q)=1$, so $r\ne p$. The result follows from Lemmas~\ref{lem:tame-hyperelliptic-reduction} and~\ref{lem:tame-hyperelliptic-irreducibility}.
\end{proof}

\subsection{Wild realisations of genus $(p^a-1)/2$}\label{ssec:wild-single-cluster}

For an integer $a\ge1$ and units $b,\eta\in\mathcal O_K^\times$, define
\begin{equation*}
f_{a,b,\eta}(x)=x^{p^a}-\pi bx-\pi\eta.
\end{equation*}
Let $X_{a,b,\eta}/K$ be the hyperelliptic curve defined by $y^2=f_{a,b,\eta}(x)$, and let $\mathcal R_{a,b,\eta}$ be the set of roots of $f_{a,b,\eta}$ in $\overline K$. When $a=1$, these curves belong to the class studied in \cite[Prop.~2.2 and Thm.~2.3]{Coppola}.

\begin{lemma}\label{lem:wild-single-reduction}
Choose a root $\alpha$ of $f_{a,b,\eta}$ and $t\in\overline K$ with $t^{2(p^a-1)}=\pi b$, and put $F=K(\alpha,t,\mu_{2(p^a-1)})$. Then the following hold.
\begin{itemize}
\item[\textup{(i)}] The extension $F/K$ is Galois and wildly ramified. The change of coordinates $x=\alpha+t^2u$ and $y=t^{p^a}w$ gives a model of $X_{a,b,\eta}$ over $F$ with good reduction and geometric special fibre
\begin{equation*}
\overline X_{a,b,\eta}:w^2=u^{p^a}-u.
\end{equation*}

\item[\textup{(ii)}] The action defined by \eqref{eq:lift-act-reduce} identifies $I(F/K)$ with
\begin{equation*}
\left\{(u,w)\mapsto(\lambda^2u+d,\lambda^{p^a}w):
d\in\F_{p^a},\ \lambda\in\mu_{2(p^a-1)}(\overline k)\right\}
\cong(C_p)^a\rtimes C_{2(p^a-1)}.
\end{equation*}

\item[\textup{(iii)}] Under this identification, $P(F/K)$ is the subgroup
\begin{equation*}
\left\{(u,w)\mapsto(u+d,w):d\in\F_{p^a}\right\}\cong(C_p)^a.
\end{equation*}
\end{itemize}
\end{lemma}

\begin{proof}
Put $f=f_{a,b,\eta}$ and $h(u)=t^{-2p^a}f(\alpha+t^2u)$. Then $h\in\mathcal O_F[u]$ and $\overline h=u^{p^a}-u$. The substitutions above give $w^2=h(u)$, which has good reduction with geometric special fibre $w^2=u^{p^a}-u$, since $p>2$ and $\overline h$ is square-free.

Put $K_0=K(\mu_{2(p^a-1)})$. The extension $K_0/K$ is unramified, and its residue field contains $\F_{p^a}$. The extensions $K_0(\alpha)/K_0$ and $K_0(t)/K_0$ are totally ramified of coprime degrees $p^a$ and $2(p^a-1)$, respectively. Their compositum $F/K_0$ is therefore totally ramified of degree $2p^a(p^a-1)$.

All roots of $\overline h$ lie in the residue field of $F$ and are simple. Hensel's lemma therefore shows that $h$ splits over $F$, and hence so does $f$. Thus $F$ is the splitting field of $f(x)(x^{2(p^a-1)}-\pi b)$ over $K$, so $F/K$ is Galois. Its ramification degree is divisible by $p$, proving \textup{(i)}.

Since $K_0/K$ is unramified and $F/K_0$ is totally ramified, we have $I(F/K)=\Gal(F/K_0)$. Moreover, $K_0(t)/K_0$ is the maximal tamely ramified subextension of $F/K_0$, so $P(F/K)=\Gal(F/K_0(t))$. For $\sigma\in I(F/K)$, the construction in \eqref{eq:lift-act-reduce} gives the automorphism
\begin{equation*}
(u,w)\mapsto(\lambda_\sigma^2u+d_\sigma,\lambda_\sigma^{p^a}w),
\end{equation*}
where $\lambda_\sigma\in\mu_{2(p^a-1)}(\overline k)$ is the reduction of $\sigma(t)/t$. The element $(\sigma(\alpha)-\alpha)/t^2$ is a root of the monic polynomial $h\in\mathcal O_F[u]$, so it is integral. Its reduction $d_\sigma$ belongs to $\F_{p^a}$, since $\overline h=u^{p^a}-u$.

This action of $I(F/K)$ on $\overline X_{a,b,\eta}$ is faithful. Indeed, an element acting trivially satisfies $\lambda_\sigma^2=\lambda_\sigma^{p^a}=1$ and $d_\sigma=0$, so $\lambda_\sigma=1$. Since reduction is injective on $\mu_{2(p^a-1)}$, this gives $\sigma(t)=t$. Moreover, $(\sigma(\alpha)-\alpha)/t^2$ is a root of $h$ reducing to $0$. Since $h(0)=0$ and $\overline h$ is square-free, Hensel's lemma shows that $0$ is the only such root, so $\sigma(\alpha)=\alpha$. As $\sigma$ also fixes $K_0$, it is the identity on $F=K_0(\alpha,t)$.

Since the action is faithful and $|I(F/K)|=2p^a(p^a-1)$, its image is the entire group displayed in \textup{(ii)}. Choose a primitive $2(p^a-1)$th root of unity $\zeta\in\overline k$. The group is generated by
\begin{equation*}
\sigma_d:(u,w)\mapsto(u+d,w),
\qquad
\tau:(u,w)\mapsto(\zeta^2u,\zeta^{p^a}w),
\end{equation*}
where $d\in\F_{p^a}$. The element $\tau$ has order $2(p^a-1)$ and satisfies $\tau\sigma_d\tau^{-1}=\sigma_{\zeta^2d}$. Thus $I(F/K)$ is the semidirect product $(C_p)^a\rtimes C_{2(p^a-1)}$. Its unique Sylow $p$-subgroup is $\{\sigma_d:d\in\F_{p^a}\}$, which therefore equals $P(F/K)$. This proves \textup{(ii)} and \textup{(iii)}.
\end{proof}

The following proposition generalises Example~\ref{ex:wild-elliptic}, which is the case $K=\Q_3$, $\pi=3$ and $a=\eta=1$.

\begin{proposition}\label{prop:wild-single-cluster}
The curve $X_{a,b,\eta}$ has genus $(p^a-1)/2$. Its Tate module $V_\ell(\Jac(X_{a,b,\eta}))$ is wildly ramified and absolutely irreducible as an $I_K$-representation, and hence as a $G_K$-representation.
\end{proposition}

\begin{proof}
Put $V=V_\ell(\Jac(X_{a,b,\eta}))\otimes_{\Q_\ell}\overline{\Q}_\ell$. Since $f_{a,b,\eta}$ is Eisenstein and therefore square-free of odd degree $p^a$, the curve has genus $(p^a-1)/2$. Let $F/K$ be the extension in Lemma~\ref{lem:wild-single-reduction}, and put $J=I(F/K)$. By part~\textup{(ii)} of that lemma, the $J$-action on $\overline X_{a,b,\eta}$ is faithful. The displayed automorphisms fix the unique point at infinity, so the identification with the action on the Tate module also applies in genus one. Remark~\ref{rem:J=I} therefore allows us to identify $J$ with $I$, and $P(F/K)$ with $P$. By Proposition~\ref{prop:inertia-automorphisms}, we may compute the action of inertia on $V$ through the action of $I$ on $\overline X_{a,b,\eta}:w^2=u^{p^a}-u$.

By Lemma~\ref{lem:wild-single-reduction}\textup{(iii)}, the elements of $P$ act by $\sigma_d:(u,w)\mapsto(u+d,w)$ for $d\in\F_{p^a}$. The quotient $\overline X_{a,b,\eta}/P$ is isomorphic to $\mathbb P^1_{\overline k}$ with parameter $w$, so $V^P=0$ by Lemma~\ref{lem:quotient-invariants}. Thus only non-trivial characters of $P$ occur in $V|_P$.

Choose a primitive $2(p^a-1)$th root of unity $\zeta\in\overline k$, and let $\tau\in I$ act by $(u,w)\mapsto(\zeta^2u,\zeta^{p^a}w)$. Under the identification $P\cong(\F_{p^a},+)$, conjugation by $\tau$ is multiplication by $\zeta^2$, which generates $\F_{p^a}^\times$. Thus $\tau$ acts on the characters of $P$ by
\begin{equation*}
(\tau\cdot\chi)(\sigma_d)
=\chi(\tau^{-1}\sigma_d\tau)
=\chi(\sigma_{\zeta^{-2}d}).
\end{equation*}
For a non-trivial character $\chi$, the equality $\tau^j\cdot\chi=\chi$ holds precisely when $\zeta^{2j}=1$. Indeed, otherwise it would imply $\chi(\sigma_{(\zeta^{-2j}-1)d})=1$ for every $d\in\F_{p^a}$, forcing $\chi$ to be trivial. Thus $I$ acts transitively on the $p^a-1$ non-trivial characters of $P$.

Since $V\ne0$ and $V^P=0$, at least one non-trivial character of $P$ occurs in $V|_P$. Since $P$ is normal in $I$, the set of characters occurring in $V|_P$ is $I$-stable, so transitivity shows that they all occur. As $\dim V=p^a-1$, each occurs with multiplicity one. Since $I$ permutes these character spaces transitively, Remark~\ref{rem:inertia-transitivity} shows that $V$ is irreducible as an $I$-representation, and hence as a $G_K$-representation. Finally, $V^P=0$ shows that $P$ acts non-trivially, so the Tate module is wildly ramified.
\end{proof}

\subsection{Wild realisations of genus $p^a-1$}\label{ssec:wild-two-clusters}

We now use the preceding family to realise the cluster picture of Corollary~\ref{cor:cluster-configurations}\textup{(ii)} and obtain curves of genus $p^a-1$. To produce two conjugate principal clusters, we take two Galois-conjugate polynomials over a quadratic extension of $K$ and translate their roots to centres with distinct reductions.

Throughout this subsection, fix $\varepsilon\in\mathcal O_K^\times$ with non-square reduction and $\delta\in\overline K$ with $\delta^2=\varepsilon$. Then $K(\delta)/K$ is unramified quadratic, and its non-trivial automorphism $\sigma$ sends $\delta$ to $-\delta$. For each integer $a\ge1$, let $D_a/K$ be the hyperelliptic curve
\begin{equation*}
D_a:y^2=f_{a,\delta,1}(x-\delta)f_{a,-\delta,1}(x+\delta).
\end{equation*}
The automorphism $\sigma$ exchanges the two factors, so their product belongs to $K[x]$. Their roots reduce to $\overline\delta$ and $-\overline\delta$, respectively, which are distinct since $p>2$. We write $\mathcal R_a$ for the set of roots of this product, which satisfies
\begin{equation*}
\mathcal R_a=(\delta+\mathcal R_{a,\delta,1})
\sqcup(-\delta+\mathcal R_{a,-\delta,1}).
\end{equation*}
Let $\alpha$ be a root of $f_{a,\delta,1}$ and $t\in\overline K$ satisfy $t^{2(p^a-1)}=\pi\delta$. We put
\begin{equation*}
F_+=K(\delta,\alpha,t,\mu_{2(p^a-1)}),
\qquad
h(u)=t^{-2p^a}f_{a,\delta,1}(\alpha+t^2u).
\end{equation*}
By Lemma~\ref{lem:wild-single-reduction}, applied over $K(\delta)$, we have $h\in\mathcal O_{F_+}[u]$ and $\overline h=u^{p^a}-u$.

\begin{lemma}\label{lem:Da-clusters}
The cluster picture of $\mathcal R_a$ is
\begin{center}
\scalebox{1.8}{\clusterpicture
  \Root{1}{first}{r1}
  \Root{}{r1}{r2}
  \RootT{}{r2}{dots1}{$\cdots$}
  \Root{}{dots1}{r3}
  \ClusterLDName s1[][\frac{1}{p^a-1}][\mathfrak s_+]
    = (r1)(r2)(dots1)(r3);
  \Root{3}{s1}{r4}
  \Root{}{r4}{r5}
  \RootT{}{r5}{dots2}{$\cdots$}
  \Root{}{dots2}{r6}
  \ClusterLDName s2[][\frac{1}{p^a-1}][\mathfrak s_-]
    = (r4)(r5)(dots2)(r6);
  \ClusterLDName R[][0][\mathcal R_a]
    = (s1)(s1n)(s2)(s2n);
\endclusterpicture}
\end{center}
Here $\mathfrak s_+=\delta+\mathcal R_{a,\delta,1}$ and $\mathfrak s_-=-\delta+\mathcal R_{a,-\delta,1}$ are the two principal clusters.
\end{lemma}

\begin{proof}
The roots of $f_{a,\delta,1}$ and $f_{a,-\delta,1}$ have valuation $1/p^a$, since both polynomials are Eisenstein over $K(\delta)$. Hence the roots in $\mathfrak s_+$ and $\mathfrak s_-$ reduce to $\overline\delta$ and $-\overline\delta$, respectively. These reductions are distinct, so $v(r-r')=0$ whenever $r\in\mathfrak s_+$ and $r'\in\mathfrak s_-$.

Since $h$ is monic and $\overline h=u^{p^a}-u$ is square-free, the roots of $h$ have distinct reductions, so their differences are units. The substitution $x=\alpha+t^2u$ multiplies these differences by $t^2$. Thus any two distinct roots of $f_{a,\delta,1}$ have difference of valuation $2v(t)=1/(p^a-1)$. The same holds for the conjugate polynomial, giving the asserted depths of $\mathfrak s_+$ and $\mathfrak s_-$.

Each $\mathfrak s_\pm$ has odd cardinality $p^a\ge3$ and only singleton children, so both are principal. The full root set $\mathcal R_a$ is not principal, since it is even and has exactly two children.
\end{proof}

Every cluster strictly contained in $\mathcal R_a$ is odd, so $\Jac(D_a)$ has potentially good reduction by \cite[Thm.~1.9(6)]{DDMM}. We may therefore apply Lemma~\ref{lem:cluster-tate}. Let $W_{\mathfrak s_\pm}$ be the representations attached to $\mathfrak s_\pm$ in that lemma. We compare their restrictions to $I_K=I_{K(\delta)}$ with the Tate modules of $X_{a,\delta,1}$ and $X_{a,-\delta,1}$.

\begin{lemma}\label{lem:Da-components}
There are isomorphisms of $I_K$-representations
\begin{equation*}
W_{\mathfrak s_+}|_{I_K}\cong V_\ell(\Jac(X_{a,\delta,1}))\otimes_{\Q_\ell}\overline{\Q}_\ell,
\quad\text{and}\quad
W_{\mathfrak s_-}|_{I_K}\cong V_\ell(\Jac(X_{a,-\delta,1}))\otimes_{\Q_\ell}\overline{\Q}_\ell.
\end{equation*}
Moreover, $W_{\mathfrak s_+}$ and $W_{\mathfrak s_-}$ are wildly ramified, irreducible and non-isomorphic as $I_K$-representations.
\end{lemma}

\begin{proof}
Choose a finite Galois extension of $K$ containing $F_+$ over which $D_a$ has semistable reduction, and let $\widetilde\Gamma_{\mathfrak s_\pm}$ denote the normalisations of the components corresponding to $\mathfrak s_\pm$. By Lemma~\ref{lem:cluster-tate}, we have
\begin{equation*}
W_{\mathfrak s_\pm}=V_\ell(\Jac(\widetilde\Gamma_{\mathfrak s_\pm}))\otimes_{\Q_\ell}\overline{\Q}_\ell.
\end{equation*}
Since $I_K$ fixes both clusters, it acts on each of these curves. The substitution $x=\delta+\alpha+t^2u'$ and $y=t^{p^a}w'$ gives
\begin{equation*}
(w')^2=h(u')f_{a,-\delta,1}(2\delta+\alpha+t^2u').
\end{equation*}
The second factor reduces to $(2\overline\delta)^{p^a}\ne0$. By \cite[Thm.~8.8]{DDMM}, the normalised component $\widetilde\Gamma_{\mathfrak s_+}$ therefore has equation $(w')^2=(2\overline\delta)^{p^a}((u')^{p^a}-u')$. Choose $c\in\overline k^\times$ with $c^2=(2\overline\delta)^{p^a}$. We obtain an isomorphism
\begin{equation*}
\phi:\widetilde\Gamma_{\mathfrak s_+}\to\overline X_{a,\delta,1},
\qquad
(u',w')\mapsto(u,w)=(u',w'/c).
\end{equation*}

Applying \eqref{eq:lift-act-reduce} to this substitution, and using that $I_K=I_{K(\delta)}$ fixes $\delta$ and acts trivially on $\overline k$, we find that $\tau\in I_K$ acts on $\widetilde\Gamma_{\mathfrak s_+}$ by
\begin{equation*}
(u',w')\mapsto(\lambda_\tau^2u'+d_\tau,\lambda_\tau^{p^a}w'),
\end{equation*}
where $\lambda_\tau=\overline{\tau(t)/t}$ and $d_\tau=\overline{(\tau(\alpha)-\alpha)/t^2}$. Since inertia fixes $c$, the isomorphism $\phi$ carries this action to $(u,w)\mapsto(\lambda_\tau^2u+d_\tau,\lambda_\tau^{p^a}w)$, which is the action on $\overline X_{a,\delta,1}$ from Lemma~\ref{lem:wild-single-reduction}. Thus $\phi$ is $I_K$-equivariant. The same argument gives an $I_K$-equivariant isomorphism $\widetilde\Gamma_{\mathfrak s_-}\cong\overline X_{a,-\delta,1}$. Proposition~\ref{prop:inertia-automorphisms}, applied over $K(\delta)$, now gives the asserted isomorphisms of representations. Proposition~\ref{prop:wild-single-cluster} then shows that both representations are wildly ramified and irreducible under $I_K$.

It remains to show that the two representations are non-isomorphic. Let $F_-/K(\delta)$ be the extension supplied by Lemma~\ref{lem:wild-single-reduction} for $X_{a,-\delta,1}$. The corresponding extension for $X_{a,\delta,1}$ is $F_+$ defined above. The actions of $I(F_\pm/K(\delta))$ on the respective Tate modules are faithful by Remark~\ref{rem:J=I}. In the genus-one case, the displayed automorphisms fix the unique point at infinity, so the same identification applies. Since $I_K=I_{K(\delta)}$, the kernels of the two $I_K$-representations are therefore $I_{F_+}$ and $I_{F_-}$.

Suppose that the representations are isomorphic. Their kernels then coincide, so their fixed fields satisfy $F_+^{\mathrm{nr}}=F_-^{\mathrm{nr}}.$ By Lemma~\ref{lem:wild-single-reduction}\textup{(ii)}, we have
\begin{equation*}
e(F_\pm/K(\delta))=2p^a(p^a-1).
\end{equation*}

Let $\beta_1,\ldots,\beta_{p^a}$ be the roots of $f_{a,-\delta,1}$. Since $\alpha\in F_+$ and the $\beta_j$ lie in $F_-$, all these roots belong to $F_+^{\mathrm{nr}}$. By Lemma~\ref{lem:Da-clusters}, we have $v(\beta_i-\beta_j)=1/(p^a-1)$ for $i\ne j$. On the other hand, $f_{a,-\delta,1}(\alpha)=2\pi\delta\alpha$, so
\begin{equation*}
\sum_{j=1}^{p^a}v(\alpha-\beta_j)
=v(f_{a,-\delta,1}(\alpha))
=1+\frac1{p^a}
<\frac{p^a}{p^a-1}.
\end{equation*}
At least one of these valuations is therefore smaller than $1/(p^a-1)$, which shows that they are all equal. This implies that $v(\alpha-\beta_j)=\frac{p^a+1}{p^{2a}}$ for every $j$. This fraction is in lowest terms. Since $\alpha-\beta_j\in F_+^{\mathrm{nr}}$ and both $F_+^{\mathrm{nr}}/F_+$ and $K(\delta)/K$ are unramified, it follows that $p^{2a}$ divides $e(F_+/K(\delta))$. This contradicts $e(F_+/K(\delta))=2p^a(p^a-1)$. This shows that the two representations are non-isomorphic.
\end{proof}

\begin{proposition}\label{prop:wild-two-clusters}
The curve $D_a$ has genus $p^a-1$, and $V_\ell(\Jac(D_a))$ is wildly ramified and absolutely irreducible.
\end{proposition}

\begin{proof}
The defining polynomial is square-free of degree $2p^a$, giving the asserted genus. We have already seen that $\Jac(D_a)$ has potentially good reduction.

Put $V=V_\ell(\Jac(D_a))\otimes_{\Q_\ell}\overline{\Q}_\ell$. Since $I_K$ fixes both principal clusters, Lemma~\ref{lem:cluster-tate}\textup{(i)} gives
\begin{equation*}
V|_{I_K}\cong W_{\mathfrak s_+}|_{I_K}\oplus W_{\mathfrak s_-}|_{I_K}.
\end{equation*}
The summands are irreducible and non-isomorphic by Lemma~\ref{lem:Da-components}. An arithmetic Frobenius lift exchanges the two clusters, since $K(\delta)/K$ is unramified quadratic. By the Galois equivariance of the cluster decomposition, it also exchanges the corresponding summands. Lemma~\ref{lem:frobenius-transitivity} therefore gives absolute irreducibility. Finally, both summands are wildly ramified by Lemma~\ref{lem:Da-components}, so the same holds for $V$.
\end{proof}

\section{Realisations by Jacobians}\label{sec:abelian-realisations}

We construct curves whose Jacobians have absolutely irreducible Tate modules of every even dimension permitted by Lemma~\ref{lem:tame}. We first state the resulting classification and its consequences. Their proofs refer to the constructions below. Throughout, $K/\Q_p$ denotes a finite extension with odd $p$.

\begin{theorem}\label{thm:classification!}
Let $g\ge1$. There exists a curve $C/K$ of genus $g$ such that $V_\ell(\Jac(C))$ is absolutely \hbox{irreducible and}
\begin{itemize}
\item[\textup{(i)}] tamely ramified if and only if $2g=\varphi(m)$ for some integer $m\ge3$ satisfying $(m,p)=1$ and \hbox{$\ord_m(q)=\varphi(m)$,}

\item[\textup{(ii)}] wildly ramified if and only if $(p-1)\mid2g$.
\end{itemize}
\end{theorem}

\begin{proof}
Necessity follows from Lemma~\ref{lem:tame}. For sufficiency in the tame case, the condition $\ord_m(q)=\varphi(m)$ implies that $(\Z/m\Z)^\times$ is cyclic. Thus $m=4$, or $m=\nu^a$ or $2\nu^a$ for an odd prime $\nu$ and $a\ge1$. In the latter two cases, $\varphi(m)=\varphi(\nu^a)$ and $\ord_m(q)=\ord_{\nu^a}(q)$. Proposition~\ref{prop:tame-superelliptic-realisation} supplies the required curve, including the case $m=4$ through $(\nu,a)=(2,2)$. For sufficiency in the wild case, apply Corollary~\ref{cor:wild-jacobian-realisations}.
\end{proof}

\begin{corollary}\label{cor:even-residue-degree}
Suppose that $K$ has even residue degree, and let $g\ge1$. There exists a curve $C/K$ of genus $g$ such that $V_\ell(\Jac(C))$ is absolutely irreducible if and only if $(p-1)\mid2g$. Moreover, every such Tate module is wildly ramified.
\end{corollary}

\begin{proof}
Since $q$ is a square, Euler's theorem gives $\ord_m(q)\mid\varphi(m)/2$ for every $m\ge3$ with $(m,p)=1$. Thus the tame condition in Theorem~\ref{thm:classification!}\textup{(i)} cannot hold. Every absolutely irreducible Tate module of a Jacobian is therefore wildly ramified in this case, and the result follows from part~\textup{(ii)}.
\end{proof}

For $X\ge1$, write $\mathcal A_K(X)$, $\mathcal J_K(X)$, $\mathcal H_K(X)$ and $\mathcal J_K$ as in \S\ref{ssec:hyper}. Comparing Theorem~\ref{thm:classification!} and \ref{thm:hyperelliptic-classification} gives the following.

\begin{corollary}\label{cor:hyperelliptic-density!}
For every $X\ge1$, we have $\mathcal J_K(X)=\mathcal A_K(X)$. Moreover, $\mathcal H_K(X)\subsetneq\mathcal J_K(X)$ for every $X\ge5(p-1)/2$.

Among the genera in $\mathcal J_K$, those admitting hyperelliptic Jacobians with absolutely irreducible Tate modules have relative density zero. More precisely,
\begin{equation*}
\lim_{X\to\infty}\frac{\#\mathcal H_K(X)}{\#\mathcal J_K(X)}=0.
\end{equation*}
Consequently, for almost every $g\in\mathcal J_K$, no Jacobian of dimension $g$ with an absolutely irreducible Tate module is $K$-isogenous to a hyperelliptic Jacobian.
\end{corollary}

\begin{proof}
The inclusion $\mathcal J_K(X)\subseteq\mathcal A_K(X)$ is immediate. Conversely, Lemmas~\ref{lem:arithmetic-consequences}\textup{(iii)} and~\ref{lem:tame} show that every $g\in\mathcal A_K(X)$ satisfies the conditions of Theorem~\ref{thm:classification!}, which provides a Jacobian of dimension $g$ with an absolutely irreducible Tate module. Thus $\mathcal J_K(X)=\mathcal A_K(X)$.

To prove the strict inclusion $\mathcal H_K(X)\subsetneq\mathcal J_K(X)$, put $g_0=5(p-1)/2$. Theorem~\ref{thm:classification!}\textup{(ii)} gives $g_0\in\mathcal J_K(X)$ whenever $X\ge g_0$. Since $(p-1)\nmid g_0$ and $p\nmid2g_0+1$, we have $g_0\ne p^a-1$ and $g_0\ne(p^a-1)/2$ for every integer $a\ge1$. This excludes the wild case. To exclude the tame case, put $r=2g_0+1=5p-4$. If $r$ is composite, there is nothing to prove. If $r$ is prime, then $5p\equiv4\pmod r$ and $r\equiv1\pmod5$, so quadratic reciprocity gives $\left(\frac{p}{r}\right)=\left(\frac{5}{r}\right)=\left(\frac{r}{5}\right)=1$. Thus $q$ is a square modulo $r$, and $\ord_r(q)\le(r-1)/2<2g_0$. Hence $g_0\notin\mathcal H_K(X)$.

For the density statement, Theorem~\ref{thm:classification!}\textup{(ii)} gives $\#\mathcal J_K(X)\ge\left\lfloor(2X)/{p-1}\right\rfloor$. The tame genera in the hyperelliptic case are counted by primes up to $2X+1$ (apart from the genus-one case), while the two wild conditions contribute $O(\log X)$ genera. By the prime number theorem, $\#\mathcal H_K(X)=O(X/\log X)$. Together with the lower bound above, this gives ${\#\mathcal H_K(X)}/{\#\mathcal J_K(X)}=O\left(1/\log X\right)$, which tends to zero.

Finally, a $K$-isogeny induces an isomorphism of rational Tate modules. Thus a Jacobian with an absolutely irreducible Tate module in a genus excluded by the hyperelliptic classification cannot be \hbox{$K$-isogenous} to a hyperelliptic Jacobian. The density result gives the final assertion.
\end{proof}

\subsection{Tame realisations by superelliptic curves}\label{ssec:tame-abelian-realisations}

Let $\nu\ne p$ be an odd prime and $a\ge1$, or let $(\nu,a)=(2,2)$. Put $d=\nu^{a-1}$, and let $\pi$ be a uniformiser of $K$. Let $Y_{\nu,a}/K$ be the superelliptic curve
\begin{equation*}
Y_{\nu,a}:y^\nu=x^{d+1}-\pi x.
\end{equation*}
This family includes the genus-one cases $(\nu,a)=(3,1)$ and $(2,2)$ used in Proposition~\ref{prop:elliptic-classification}. The following lemma gives good reduction after a tame extension and computes the inertia action using \S\ref{ssec:inertia-automorphisms}.

\begin{lemma}\label{lem:tame-superelliptic-reduction}
The following hold.
\begin{itemize}
\item[\textup{(i)}] The curve $Y_{\nu,a}$ has genus $\varphi(\nu^a)/2$.

\item[\textup{(ii)}] Choose $t\in\overline K$ with $t^{\nu^a}=\pi$, let $\zeta$ be a primitive $\nu^a$th root of unity, and put $F=K(t,\zeta)$. Then $F/K$ is Galois and tamely ramified, with $I(F/K)\cong C_{\nu^a}$. The substitution $x=t^\nu u$ and $y=t^{d+1}w$ gives
\begin{equation*}
w^\nu=u^{d+1}-u,
\end{equation*}
which has good reduction over $F$.

\item[\textup{(iii)}] Let $\overline Y_{\nu,a}/\overline k$ be the geometric special fibre, and let $\tau$ be the generator of $I(F/K)$ satisfying $\tau(t)=\zeta t$. Under \eqref{eq:lift-act-reduce}, $\tau$ acts on $\overline Y_{\nu,a}$ via
\begin{equation*}
(u,w)\mapsto(\bar\zeta^\nu u,\bar\zeta^{d+1}w),
\end{equation*}
where $\bar\zeta$ denotes the reduction of $\zeta$ modulo $\mathfrak m_F$.
\end{itemize}
\end{lemma}
\begin{proof}
The polynomial $x^{d+1}-\pi x=x(x^d-\pi)$ is square-free of degree $d+1$. Since $\gcd(\nu,d+1)=1$, the genus formula for superelliptic curves \cite[Prop.~2(3)]{GPS} gives $g(Y_{\nu,a})=\frac{(\nu-1)d}{2}=\frac{\varphi(\nu^a)}{2}$, proving \textup{(i)}.

Since $p\nmid\nu^a$, the extension $K(\zeta)/K$ is unramified. The polynomial $X^{\nu^a}-\pi$ is Eisenstein over $K(\zeta)$, so $F/K(\zeta)$ is totally ramified of degree $\nu^a$. Moreover, $F$ is the splitting field of this polynomial over $K$. Thus $F/K$ is Galois and tamely ramified, and $I(F/K)=\Gal(F/K(\zeta))$ is cyclic of order $\nu^a$, generated by $t\mapsto\zeta t$.

The substitution $x=t^\nu u$ and $y=t^{d+1}w$ gives $w^\nu=u^{d+1}-u$. For $a\ge2$, the polynomial $u^{d+1}-u$ has discriminant $\pm d^d$. Since $p\nmid\nu d$, this equation has good reduction over $F$ by \cite[\S1]{Sutherland}.

When $a=1$, the equivalent equation $(2u-1)^2=1+4w^\nu$ gives the same conclusion. Indeed, the polynomial $1+4w^\nu$ has degree $\nu\ge3$, unit leading coefficient and square-free reduction, since $p\nmid2\nu$. Thus \cite[\S1]{Sutherland} gives good reduction. This proves \textup{(ii)}.

Part~\textup{(iii)} follows by applying \eqref{eq:lift-act-reduce} to this substitution, exactly as in Example~\ref{ex:inertia_action_tame}.
\end{proof}

Fix a prime $\ell\ne p$ and put $V=V_\ell(\Jac(Y_{\nu,a}))\otimes_{\Q_\ell}\overline{\Q}_\ell$.

\begin{lemma}\label{lem:tame-superelliptic-characters}
The action of $I_K$ on $V$ factors through a faithful action of $I=I(F/K)\cong C_{\nu^a}$, and
\begin{equation*}
V|_I\cong\bigoplus_{\chi\text{ faithful}}\chi,
\end{equation*}
where the sum runs over all faithful characters of $I$, each occurring with multiplicity one.
\end{lemma}

\begin{proof}
By Lemma~\ref{lem:tame-superelliptic-reduction}, the curve $Y_{\nu,a}$ has good reduction over $F$, so the action of $I_K$ on $V$ factors through $I=I(F/K)\cong C_{\nu^a}$. Let $\tau$ be the generator of $I$ satisfying $\tau(t)=\zeta t$, and put $H=\langle\tau^d\rangle$, the subgroup of order $\nu$. By Lemma~\ref{lem:tame-superelliptic-reduction}\textup{(iii)}, $H$ fixes $u$ and multiplies $w$ by the $\nu$th roots of unity. Hence $\overline Y_{\nu,a}/H\cong\mathbb P^1_{\overline k}$ with parameter $u$. Proposition~\ref{prop:inertia-automorphisms} and Lemma~\ref{lem:quotient-invariants}\textup{(ii)} therefore give $V^H=0$.

Every non-faithful character of the cyclic group $I$ is trivial on $H$, since every non-trivial subgroup of $I$ contains $H$. Thus only faithful characters of $I$ occur in $V|_I$. By Lemma~\ref{lem:quotient-invariants}\textup{(i)}, the trace of every element of $I$ on $V$ is rational. The group $\Gal(\Q(\mu_{\nu^a})/\Q)$ permutes the faithful characters of $I$ transitively, so their multiplicities in $V|_I$ are equal. Finally, $\dim V=\varphi(\nu^a)$ by Lemma~\ref{lem:tame-superelliptic-reduction}\textup{(i)}, and $I=C_{\nu^a}$ has exactly $\varphi(\nu^a)$ faithful characters. Each therefore occurs with multiplicity one. In particular, the action of $I$ on $V$ is faithful.
\end{proof}

\begin{proposition}\label{prop:tame-superelliptic-realisation}
The curve $Y_{\nu,a}$ has genus $\varphi(\nu^a)/2$. Its Tate module $V_\ell(\Jac(Y_{\nu,a}))$ is tamely ramified and is absolutely irreducible if and only if $\ord_{\nu^a}(q)=\varphi(\nu^a)$.
\end{proposition}

\begin{proof}
The genus and tame ramification assertions follow from Lemma~\ref{lem:tame-superelliptic-reduction}. By Lemma~\ref{lem:tame-superelliptic-characters}, the restriction $V|_I$ contains each faithful character of $I$ exactly once and no other characters.

An arithmetic Frobenius lift $\Phi$ satisfies $\Phi\tau\Phi^{-1}=\tau^q$ and therefore permutes these characters by $\chi\mapsto\chi^{q^{-1}}$, where $q^{-1}$ is taken modulo $\nu^a$. Every Frobenius orbit therefore has length $\ord_{\nu^a}(q)$. By Lemma~\ref{lem:frobenius-transitivity}, $V$ is irreducible if and only if $\ord_{\nu^a}(q)=\varphi(\nu^a)$.
\end{proof}

\subsection{Wild realisations by Weil restriction}\label{ssec:wild-weil-realisations}

Taking $a=1$ in the family of \S\ref{ssec:wild-single-cluster} gives hyperelliptic Jacobians of dimension $(p-1)/2$ whose Tate modules are wildly ramified and absolutely irreducible as $I_K$-representations. To obtain every positive multiple of this dimension, we use Weil restriction from unramified extensions of $K$.

Fix $r\ge1$, and let $K_r/K$ be the unramified extension of degree $r$, with residue field $k_r$. The uniformiser $\pi$ of $K$ remains a uniformiser of $K_r$. Choose $b\in\mathcal O_{K_r}^\times$ such that $k_r=k(\bar b)$, and consider
\begin{equation*}
X_b=X_{1,b,1}:y^2=x^p-\pi bx-\pi
\end{equation*}
over $K_r$. Set $A_r=\Res_{K_r/K}\Jac(X_b)$.

\begin{lemma}\label{lem:wild-distinct-inertia}
Let $b,b'\in\mathcal O_{K_r}^\times$ satisfy $\bar b\ne\bar b'$. Then the representations
\begin{equation*}
V_\ell(\Jac(X_b))\otimes_{\Q_\ell}\overline{\Q}_\ell
\qquad\text{and}\qquad
V_\ell(\Jac(X_{b'}))\otimes_{\Q_\ell}\overline{\Q}_\ell
\end{equation*}
are irreducible and non-isomorphic as $I_K$-representations.
\end{lemma}

\begin{proof}
Their irreducibility follows from Proposition~\ref{prop:wild-single-cluster}, applied over $K_r$. Since $K_r/K$ is unramified, we have $I_{K_r}=I_K$ and $P_{K_r}=P_K$. Write $V_b=V_\ell(\Jac(X_b))\otimes_{\Q_\ell}\overline{\Q}_\ell$, and let $\mathcal R_b$ denote the roots of $x^p-\pi bx-\pi$. Define $V_{b'}$ and $\mathcal R_{b'}$ similarly.

It remains to show that $V_b$ and $V_{b'}$ are non-isomorphic. We apply the argument of Lemma~\ref{lem:Da-components} over $K_r$, with $a=1$ and $\delta,-\delta$ replaced by $b,b'$. Let $F_b/K_r$ and $F_{b'}/K_r$ be the extensions supplied by Lemma~\ref{lem:wild-single-reduction} for $X_b$ and $X_{b'}$, respectively. As before, an isomorphism $V_b\cong V_{b'}$ would give $F_b^{\mathrm{nr}}=F_{b'}^{\mathrm{nr}}$, where $e(F_b/K_r)=2p(p-1)$. However, since $b-b'$ is a unit, the same valuation calculation gives $v(\alpha-\beta)=(p+1)/p^2$ for roots $\alpha$ of $f_{1,b,1}$ and $\beta$ of $f_{1,b',1}$. Both roots lie in $F_b^{\mathrm{nr}}$, forcing $p^2\mid e(F_b/K_r)$, a contradiction.
\end{proof}

\begin{proposition}\label{prop:wild-weil-realisation}
Let $r\ge1$, let $K_r/K$ be the unramified extension of degree $r$, and choose $b\in\mathcal O_{K_r}^\times$ such that $k_r=k(\bar b)$. Then $A_r$ has dimension $r(p-1)/2$, and $V_\ell(A_r)$ is wildly ramified and absolutely irreducible.
\end{proposition}

\begin{proof}
By Proposition~\ref{prop:wild-single-cluster}, the curve $X_b$ has genus $(p-1)/2$, so $\dim A_r=r(p-1)/2$.

Let $V=V_\ell(A_r)\otimes\overline{\Q}_\ell$, and let $\Phi\in G_K$ be an arithmetic Frobenius lift. Since $K_r/K$ is unramified of degree $r$, its restriction $\sigma=\Phi|_{K_r}$ generates $\Gal(K_r/K)$. For $0\le i<r$, put $b_i=\sigma^i(b)$ and \mbox{$V_i=V_\ell(\Jac(X_{b_i}))\otimes_{\Q_\ell}\overline{\Q}_\ell$.} The Weil-restriction decomposition \cite[\S1.2, (1)]{DiemNaumann} gives
\begin{equation*}
(A_r)_{K_r}\cong\prod_{i=0}^{r-1}\Jac(X_{b_i}),
\qquad
V|_{G_{K_r}}\cong\bigoplus_{i=0}^{r-1}V_i.
\end{equation*}

Under this decomposition, points on $A_r(\overline K)$ correspond to tuples $(P_0,\ldots,P_{r-1})$ with $P_i\in\Jac(X_{b_i})(\overline K)$. If $h\in G_K$ satisfies $h|_{K_r}=\sigma^j$, then $h(b_i)=b_{i+j}$ and $h(\pi)=\pi$. Therefore, applying $h$ to a point $(x,y)\in X_{b_i}(\overline K)$ gives a point $(h(x),h(y))\in X_{b_{i+j}}(\overline K)$. Passing to Jacobians gives a map $P_i\mapsto hP_i$ from the $i$th factor to the $(i+j)$th factor. Thus the action of $G_K$ on the product is given by
\begin{equation*}
h(P_0,\ldots,P_{r-1})=(hP_{-j},hP_{1-j},\ldots,hP_{r-1-j}),
\end{equation*}
with indices taken modulo $r$; see \cite[\S1.2, (1)--(2)]{DiemNaumann}.

For $h\in G_{K_r}$, we have $j=0$, so $h$ acts coordinatewise and preserves each factor. Thus each $V_i$ is $G_{K_r}$-stable and hence $I_K$-stable, since $I_K=I_{K_r}\subseteq G_{K_r}$. On the other hand, $\Phi|_{K_r}=\sigma$, so $\Phi$ sends the $i$th factor to the $(i+1)$st. Passing to Tate modules gives $\Phi(V_i)=V_{i+1}$, with indices modulo $r$. Thus $\Phi$ permutes the $V_i$ transitively.

By Proposition~\ref{prop:wild-single-cluster}, each $V_i$ is irreducible as an $I_K$-representation. Moreover, $\bar b_i=\bar b^{q^i}$ are distinct, since $k_r=k(\bar b)$ and $[k_r:k]=r$. Lemma~\ref{lem:wild-distinct-inertia} therefore shows that the $V_i$ are pairwise non-isomorphic as $I_K$-representations. Lemma~\ref{lem:frobenius-transitivity} now shows that $V$ is irreducible. Finally, Proposition~\ref{prop:wild-single-cluster} shows that $P_K$ acts non-trivially on each $V_i$, so $V_\ell(A_r)$ is wildly ramified.
\end{proof}

To pass from a Weil restriction to a Jacobian, we use the construction suggested in \cite{MO515765}. The following lemma allows arbitrary residue characteristic, as it will be used over finite extensions $L/\Q_2$ in \S\ref{sec:dyadic}.

\begin{lemma}\label{lem:compact-type-smoothing}
Let $K$ be a finite extension of $\Q_p$, where $p$ is any prime, and let $M/K$ be a finite extension. Let $X/M$ be a curve of genus $h\ge1$ with an $M$-rational point. There exists a curve $C/K$ of genus $[M:K]h$ such that
\begin{equation*}
V_\ell(\Jac(C))\cong V_\ell(\Res_{M/K}\Jac(X))
\end{equation*}
as $G_K$-representations for every prime $\ell\ne p$.
\end{lemma}

\begin{proof}
By \cite[Cor.~1.6 and proof of Thm.~5.1]{DDWeil}, it suffices to prove this assertion for a fixed prime $\ell\ne p$. Put $r=[M:K]$ and $A=\Res_{M/K}\Jac(X)$. If $r=1$, take $C=X$, so suppose that $r\ge2$.

Choose $P\in X(M)$ and a closed point $Q\in\mathbb P^1_K$ with residue field $M$. Regard $X$ as a curve over $K$ via $X\to\operatorname{Spec}M\to\operatorname{Spec}K$, and glue it to $\mathbb P^1_K$ by identifying $P$ with $Q$. Over $\overline K$, this joins the projective line to each of the $r$ conjugates of $X$ at the corresponding conjugate of $P$. If $r\ge3$, the resulting curve is stable, meaning that it is nodal and has only finitely many automorphisms over $\overline K$. Indeed, the projective line has at least three attachment points, and each conjugate of $X$ has one; see \cite[Tag~0E73]{Stacks}. If $r=2$, the projective line has only two attachment points, so it admits infinitely many automorphisms fixing them. We therefore contract this component as in \cite[Tag~0E3M]{Stacks}. In either case, write $\Gamma/K$ for the resulting stable curve. The geometric components of $\Gamma$ meet without forming cycles, so $\Gamma$ is of compact type. Consequently, over $\overline K$, $\textup{Pic}^0\Gamma$ is the product of the Jacobians of the conjugates of $X$, since the projective line does not contribute. This gives $\operatorname{Pic}^0(\Gamma)\cong A$ over $K$.

By \cite[Cor.~1.7 and Cor.~1.9]{DeligneMumford}, there is a projective family $\mathcal X\to S$ of stable curves, where $S$ is a smooth irreducible $K$-variety, such that $\mathcal X_{s_0}\cong\Gamma$ for some $s_0\in S(K)$ and the fibres are smooth outside a proper closed subset of $S$. After shrinking $S$ around $s_0$, we may assume that every fibre is of compact type. Their Jacobians then form an abelian scheme $\mathcal J\to S$ with $\mathcal J_{s_0}\cong A$; see \cite[\S\S1.2--1.3]{MoonenOort}.

By \cite[Thm. 5.1]{KisinLocalConstancy}, applied to $\mathcal J\to S$, there is a neighbourhood $U$ of $s_0$ in $S(K)$ with $V_\ell(\mathcal J_s)\cong V_\ell(A)$ as $G_K$-representations for every $s\in U$. Since $S$ is smooth and the non-smooth fibres occur only over a proper closed subset, we can choose $s\in U$ for which $\mathcal X_s$ is smooth. Then $C=\mathcal X_s$ is a smooth projective geometrically connected curve satisfying $V_\ell(\Jac(C))\cong V_\ell(A)$, as required.
\end{proof}

\begin{corollary}\label{cor:wild-jacobian-realisations}
Let $g\ge1$. There exists a curve $C/K$ of genus $g$ such that $V_\ell(\Jac(C))$ is wildly ramified and absolutely irreducible if and only if $(p-1)\mid2g$.
\end{corollary}

\begin{proof}
Necessity follows from Lemma~\ref{lem:tame}\textup{(ii)}. Conversely, put $r=2g/(p-1)$ and choose $X_b/K_r$ as in Proposition~\ref{prop:wild-weil-realisation}. Since $X_b$ has a $K_r$-rational point at infinity, Lemma~\ref{lem:compact-type-smoothing} gives a curve $C/K$ of genus $g$ with $V_\ell(\Jac(C))\cong V_\ell(A_r)$. The latter is wildly ramified and absolutely irreducible by Proposition~\ref{prop:wild-weil-realisation}.
\end{proof}

\section{Genus gap for curves with potentially good reduction}
\label{sec:genus-gap}

In the tame case, the constructions of \S\ref{ssec:tame-abelian-realisations} give curves with potentially good reduction whose Tate modules account for every dimension in Theorem~\ref{thm:classification}\textup{(i)}. In the wild case of Theorem~\ref{thm:classification}, we ask which dimensions can be realised by $V_\ell(\Jac(C))$ when the curve $C$ itself is required to have potentially good reduction. When $K$ has even residue degree and odd residue characteristic, we show that no genus strictly between $(p-1)/2$ and $p(p-1)/2$ can occur. We first record the structure of the centraliser of wild inertia. We write $D=C_I(P)$ for this group.

\begin{lemma}\label{lem:abelian-inertia-structure}
Suppose that $V$ is absolutely irreducible and that $P$ is abelian. Then $D\cong P\times T$, where $T$ is cyclic of order prime to $p$ and central in $I$. Moreover, $D$ and $T$ are normal in $G$.
\end{lemma}

\begin{proof}
By Schur--Zassenhaus, we may write $I=P\rtimes J$, where $J\cong I/P$ is cyclic of order prime to $p$. Since $P$ is abelian, we have $D=P\times C_J(P)$. Put $T=C_J(P)$. This subgroup is cyclic and commutes with both $P$ and $J$, so it is central in $I$. Since $I$ and $P$ are normal in $G$, we have $\sigma D\sigma^{-1}=C_{\sigma I\sigma^{-1}}(\sigma P\sigma^{-1})=C_I(P)=D$ for every $\sigma\in G$. Thus $D$ is normal in $G$. Finally, $T$ consists of the elements of $D$ of order prime to $p$, so it is characteristic in $D$ and hence normal in $G$.
\end{proof}

\begin{theorem}\label{thm:gap!}
Let $K/\Q_p$ be finite with $p>2$ and even residue degree. Let $C/K$ be a smooth projective geometrically connected curve of genus $g\ge1$ with potentially good reduction. If $V_\ell(\Jac(C))$ is absolutely irreducible, then
\begin{equation*}
g=\frac{p-1}2\qquad\text{or}\qquad g\ge\frac{p(p-1)}2.
\end{equation*}
\end{theorem}

\begin{proof}
By Corollary~\ref{cor:even-residue-degree}, $V$ is wildly ramified. If $g=1$, Lemma~\ref{lem:tame}\textup{(ii)} gives $(p-1)\mid2$, so $p=3$ and $g=(p-1)/2$.

Suppose that $2\le g<p(p-1)/2$. Choose a finite Galois extension $F/K$ over which $C$ has good reduction. By Proposition~\ref{prop:inertia-automorphisms}, $P$ acts faithfully on the geometric special fibre $\overline C$. If $|P|\ge p^2$, then $P$ contains a subgroup $H$ of order $p^2$, which is abelian. By \cite[Cor.~2.9]{MatignonRocher}, $p^2=|H|\le\frac{2pg}{p-1}.$ This violates the assumed bound on $g$. Since $P$ is non-trivial, it follows that $P\cong C_p$.

Write $D=C_I(P)=P\times T$ as in Lemma~\ref{lem:abelian-inertia-structure}, and put $s=|T|$. Since $P$ is normal in $G$ and acts non-trivially on $V$, irreducibility gives $V^P=0$. Thus every character occurring in $V|_D$ has the form $\theta\otimes\psi$, where $\theta$ is a non-trivial character of $P$ and $\psi$ is a character of $T$. As $\Aut(P)\cong\F_p^\times$ acts freely on the non-trivial characters of $P$, the stabiliser of $\theta\otimes\psi$ in $I$ is $C_I(P)=D$. Hence $\Ind_D^I(\theta\otimes\psi)$ is irreducible, and Frobenius reciprocity and Lemma~\ref{lem:multiplicity-free} show that $\theta\otimes\psi$ occurs in $V|_D$ with multiplicity at most one.

Let $\Phi$ be an arithmetic Frobenius lift. Since $T$ is normal in $G$ and embeds into $I/P$, we have $\Phi t\Phi^{-1}=t^q$ for every $t\in T$. As $T$ is central in $I$, it acts on each irreducible constituent of $V|_I$ through a character. By Lemma~\ref{lem:multiplicity-free}, $\Phi$ permutes these constituents transitively, so the characters of $T$ occurring in $V$ form a single orbit under $\psi\mapsto\psi^q$. In particular, they have a common kernel. This kernel acts trivially on $V$, so these characters are faithful. Since $V|_T$ has rational character by Lemma~\ref{lem:quotient-invariants}, all $\varphi(s)$ faithful characters of $T$ occur. Hence $\ord_s(q)=\varphi(s)$.

Since $K$ has even residue degree, $q$ is a square. As a square cannot generate $(\Z/s\Z)^\times$ when $s\ge3$, we obtain $s\le2$. Thus $T$ has a unique faithful character $\psi_0$, and every character occurring in $V|_D$ has the form $\theta\otimes\psi_0$ with $\theta\ne\mathds{1}$, with multiplicity at most one. Hence $2g=\dim V\le p-1$. Since $(p-1)\mid2g$ by Lemma~\ref{lem:tame}\textup{(ii)}, we conclude that $g=(p-1)/2$.
\end{proof}

The genus gap gives isogeny classes of Jacobians with potentially good reduction that contain no Jacobian of a curve with potentially good reduction.

\begin{corollary}\label{cor:no-potential-good-isogeny!}
Suppose that $K$ has even residue degree. For every integer $r$ with $2\le r\le p-1$, put $g=r(p-1)/2$. There exists a curve $C/K$ of genus $g$ with absolutely irreducible Tate module such that $\Jac(C)$ has potentially good reduction, but its $K$-isogeny class contains no Jacobian of a curve with potentially good reduction.
\end{corollary}

\begin{proof}
Theorem~\ref{thm:classification} gives such a curve $C$, and Lemma~\ref{lem:arithmetic-consequences} gives potentially good reduction of $\Jac(C)$. Every abelian variety that is $K$-isogenous to $\Jac(C)$ has the same Tate module and dimension. Since ${p-1}/{2}<g<{p(p-1)}/{2},$ Theorem~\ref{thm:gap!} rules out potentially good reduction for every curve whose Jacobian belongs to this isogeny class.
\end{proof}

We conclude this section by noting that both endpoints of the gap in Theorem~\ref{thm:gap!} are attained. By Lemma~\ref{lem:wild-single-reduction} and Proposition~\ref{prop:wild-single-cluster}, the curves $X_{1,b,\eta}$ realise the lower endpoint $(p-1)/2$ over every finite extension $K/\Q_p$. The following lemma gives superelliptic curves attaining the upper endpoint $p(p-1)/2$ over every finite unramified extension of $\Q_p$, with no assumption on the residue degree.

\begin{lemma}\label{lem:gap-sharpness}
Let $K/\Q_p$ be a finite unramified extension. The curve $C/K$ defined by $y^p=x^{p+1}+x^p+1$ has genus $p(p-1)/2$ and potentially good reduction. Its Tate module $V_\ell(\Jac(C))$ is absolutely irreducible.
\end{lemma}

\begin{proof}
The polynomial $x^{p+1}+x^p+1$ is square-free of degree $p+1$, so the genus formula for superelliptic curves \cite[Prop.~2(3)]{GPS} gives $g(C)=p(p-1)/2$. Put $V=V_\ell(\Jac(C))\otimes_{\Q_\ell}\overline{\Q}_\ell$. The hypotheses of \cite[Thm.~3.1]{ChretienMatignon} hold for $n=1$ and $c=1$, since $c$ is a unit with $c^p=c$. That theorem and Step~E of its proof show that $C$ has potentially good reduction, that $V^{P_K}=0$, and that $\operatorname{Sw}(V)=1$ over $\Q_p^{\mathrm{nr}}$. The auxiliary extension used in their calculation is tame, so it does not change wild inertia. Since $K/\Q_p$ is unramified, $K$ and $\Q_p$ have the same inertia group, so $\operatorname{Sw}(V)=1$ over $K$ as well.

Since $V^{P_K}=0$, every irreducible constituent of $V$ has positive integral Swan conductor, so additivity and $\operatorname{Sw}(V)=1$ force $V$ to be irreducible.
\end{proof}

\section{Curves in residue characteristic two}\label{ssec:dyadic-realisations}\label{sec:dyadic}

We now complete the proof of Theorem~\ref{thm:classification} in residue characteristic two and obtain the representation-theoretic classification in Corollary~\ref{cor:abstract}.

Throughout this section, let $L/\Q_2$ be a finite extension with residue field $k_L$ of cardinality $q_L$. Choose a uniformiser $\pi_L$, and write $v_L$ for the valuation normalised by $v_L(\pi_L)=1$.

\begin{theorem}\label{thm:dyadic-jacobian-realisations}
Let $g\ge1$ and let $\ell$ be an odd prime. Then the following hold.
\begin{itemize}
\item[\textup{(i)}] There exists a curve $C/L$ of genus $g$ such that $V_\ell(\Jac(C))$ is tamely ramified and absolutely irreducible if and only if $2g=\varphi(m)=\ord_m(q_L)$ for some odd integer $m\ge3$.

\item[\textup{(ii)}] There exists a curve $C/L$ of genus $g$ such that $V_\ell(\Jac(C))$ is wildly ramified and absolutely irreducible.
\end{itemize}
\end{theorem}
Combining the constructions below with Theorem~\ref{thm:classification!} and Lemma~\ref{lem:tame} gives the following classification of the possible dimensions of absolutely irreducible Galois representations for which the trace of every element of inertia is rational.
\begin{corollary}\label{cor:abstract}
Let $\mathcal K/\Q_p$ be a finite extension with residue field of cardinality $q$, and let $\ell\ne p$ be prime. The dimensions of absolutely irreducible $\ell$-adic representations of $G_{\mathcal K}$ for which the trace of every element of inertia is rational are precisely
\begin{equation*}
\{1\}\cup\{n\ge1:(p-1)\mid n\}\cup
\{\varphi(m):m\ge3,\ (m,p)=1,\ \ord_m(q)=\varphi(m)\}.
\end{equation*}
Every even dimension in this set occurs as the dimension of an absolutely irreducible Tate module $V_\ell(\Jac(C))$ for a curve $C/\mathcal K$. When $\mathcal K$ has even residue degree, the last set in the union is empty.
\end{corollary}

\begin{proof}
Necessity follows from Lemma~\ref{lem:tame}. For $p=2$, the condition $(p-1)\mid n$ is vacuous, and every positive dimension is realised in \S\ref{ssec:dyadic-representations}. The assertion concerning Jacobians follows from Theorem~\ref{thm:dyadic-jacobian-realisations}\textup{(ii)}.

For $p>2$, dimension $1$ is realised by the trivial representation. Every other dimension in the set is even and is realised by a Jacobian by Theorem~\ref{thm:classification!}. The final assertion follows from Lemma~\ref{lem:tame}, since $q$ is then a square.
\end{proof}

Part~\textup{(i)} of the main theorem of this section follows from the constructions of \S\ref{ssec:tame-abelian-realisations}, with a minor adjustment to the proof.

\begin{proof}[Proof of Theorem~\ref{thm:dyadic-jacobian-realisations}\textup{(i)}]
Necessity follows from Lemmas~\ref{lem:arithmetic-consequences} and~\ref{lem:tame}. Conversely, suppose that $m\ge3$ is odd and $2g=\varphi(m)=\textup{ord}_m(q_L)$. Since $(\Z/m\Z)^\times$ is cyclic, we have $m=\nu^a$ for an odd prime $\nu$ and an integer $a\ge1$. Consider
\begin{equation*}
Y_{\nu,a}:y^\nu=x^{\nu^{a-1}+1}-\pi_Lx.
\end{equation*}
The arguments of \S\ref{ssec:tame-abelian-realisations} apply verbatim over $L$, except for the proof of good reduction when $a=1$. That proof uses the model $z^2=1+4w^\nu$, obtained by putting $z=2u-1$. Its reduction in characteristic two is $(z-1)^2=0$, which is not smooth. We can instead choose $t$ with $t^\nu=\pi_L$ and put $F=L(t,\mu_\nu)$. The substitution $x=t^\nu u$, $y=t^2w$ gives the equation $u^2-u=w^\nu$ of a hyperelliptic curve whose discriminant $\pm\nu^\nu$ is a unit. As a result, the geometric special fibre is $\overline Y_{\nu,1}:u^2+u=w^\nu$.

The subsequent computations of inertia and Frobenius are unchanged. Thus $Y_{\nu,a}$ has genus $g$ and a tamely ramified, absolutely irreducible Tate module for every prime $\ell\ne2$.
\end{proof}

We turn to the wild case in the following subsections, first constructing representations of every positive dimension and then realising every positive genus by a Jacobian.

\subsection{Representations of every dimension}\label{ssec:dyadic-representations}

For $r\ge1$, let $L_r/L$ be the unramified extension of degree $r$, with residue field $k_{L_r}$. The element $\pi_L$ remains a uniformiser of $L_r$, and we have $I_{L_r}=I_L$ and $P_{L_r}=P_L$. We also write $v_L$ for its extension to $\overline L$. We first construct quadratic characters whose restrictions to inertia are distinct. These give the representations needed for Corollary~\ref{cor:abstract} and the construction in \S\ref{ssec:dyadic-jacobians}.

\begin{lemma}\label{lem:dyadic-quadratic-characters}
For $b\in\mathcal O_{L_r}^\times$, let $\eta_b:G_{L_r}\to\{\pm1\}$ be the quadratic character of $L_r(\sqrt{1+\pi_Lb})/L_r$. Then $\eta_b$ is ramified. If $\bar b\ne\bar b'$, then $\eta_b|_{I_L}\ne\eta_{b'}|_{I_L}$.
\end{lemma}

\begin{proof}
Let $u\in\mathcal O_{L_r}^\times$ satisfy $v_L(u-1)=1$. We first show that $u$ is not a square in any finite unramified extension of $L_r$. If $u=z^2$ in such an extension, then $\bar z^2=1$ in characteristic two, so $\bar z=1$. Writing $z=1+c$ with $v_L(c)\ge1$ gives $v_L(u-1)=v_L(2c+c^2)\ge2$, a contradiction.

Applying this to $u=1+\pi_Lb$ shows that $\eta_b$ is ramified. If $\bar b\ne\bar b'$, then $\frac{1+\pi_Lb}{1+\pi_Lb'}-1=\frac{\pi_L(b-b')}{1+\pi_Lb'}$ has valuation one. The same argument shows that $\eta_b\eta_{b'}^{-1}$ is ramified, so their restrictions to inertia are distinct.
\end{proof}

\begin{proposition}\label{prop:dyadic-representations}
For every prime $\ell\ne2$ and every integer $r\ge1$, there exists a wildly ramified, absolutely irreducible continuous $\Q_\ell$-representation of $G_L$ of dimension $r$ for which the trace of every element of inertia is rational.
\end{proposition}

\begin{proof}
Choose $b\in\mathcal O_{L_r}^\times$ such that $k_{L_r}=k_L(\bar b)$, and put $W=\Ind_{G_{L_r}}^{G_L}\eta_b$. Then $\dim W=r$. Let $\Phi\in G_L$ be an arithmetic Frobenius lift, and put $b_i=\Phi^i(b)$ for $0\le i<r$. We have
\begin{equation*}
W|_{I_L}\cong\bigoplus_{i=0}^{r-1}\eta_{b_i}|_{I_L}.
\end{equation*}
Since $\bar b_i=\bar b^{q_L^i}$ are distinct, Lemma~\ref{lem:dyadic-quadratic-characters} shows that these characters are pairwise non-isomorphic. Frobenius permutes them cyclically, and Lemma~\ref{lem:frobenius-transitivity} gives absolute irreducibility.

For every $h\in I_L$, its trace on $W$ is $\sum_{i=0}^{r-1}\eta_{b_i}(h)\in\Z$. Finally, each $\eta_{b_i}$ is ramified and quadratic, hence wildly ramified in residue characteristic two. Thus $W$ is wildly ramified.
\end{proof}

\subsection{Wild realisations by Jacobians}\label{ssec:dyadic-jacobians}

We construct abelian varieties of every positive dimension over $L$ by taking Weil restrictions of suitable quadratic twists of an elliptic curve. Lemma~\ref{lem:compact-type-smoothing} then gives curves whose Jacobians have the same Tate modules.

For an elliptic curve over $L$ with potentially good reduction, the non-trivial possibilities for the inertia image are
\begin{equation*}
C_2,\ C_3,\ C_4,\ C_6,\ Q_8,\ \mathrm{SL}_2(\F_3)
\end{equation*}
by \cite[\S 1, pp. 354--355]{Kraus}. In either non-abelian case, the wild inertia image is $Q_8$. We use an elliptic curve with inertia image $\mathrm{SL}_2(\F_3)$. Twisting by the quadratic characters $\eta_b$ constructed in \S\ref{ssec:dyadic-representations} gives elliptic curves over $L_r$ whose Tate modules are irreducible under inertia and pairwise non-isomorphic for distinct reductions of $b$.

\begin{lemma}\label{lem:dyadic-elliptic}
There exists an elliptic curve $E/L$ with potentially good reduction and \hbox{inertia image} \hbox{$I\cong\mathrm{SL}_2(\F_3)$.} For every odd prime $\ell$, $V_\ell(E)$ is wildly ramified and absolutely irreducible as a $P_L$-representation, and hence as an $I_L$-representation and a $G_L$-representation.
\end{lemma}

\begin{proof}
Consider $E:y^2=x^3+x^2+64/\pi_L$, and write $c_4,c_6,\Delta$ and $j$ for its usual invariants. A direct computation gives $c_4=16$, $v_L(\Delta)=12v_L(2)-1$ and $v_L(j)=1$. Thus $E$ has potentially good reduction. Since changing a Weierstrass equation changes the valuation of its discriminant by a multiple of twelve, $E$ cannot have good reduction over $L$. Its reduction is therefore additive.

Put $\gamma(X)=X^8-6c_4X^4-8c_6X^2-3c_4^2$. The polynomial $\pi_LT^8\gamma(2+4/T)/2^{16}$ is Eisenstein at $\pi_L$, so $\gamma$ is irreducible over $L$. Since $c_4\ne0$ and $3\nmid v_L(\Delta)$, \cite[Prop.~2 and Lem.~6]{DDRootNumbers} gives $I\cong\mathrm{SL}_2(\F_3)$.

Fix an odd prime $\ell$ and put $V=V_\ell(E)\otimes_{\Q_\ell}\overline{\Q}_\ell$. The wild inertia image is the Sylow $2$-subgroup $P\cong Q_8$ of $I$. Its irreducible representations consist of four one-dimensional characters $\mathds{1},\chi_2,\chi_3,\chi_4$ and one two-dimensional representation $\sigma$.

Conjugation by $I/P\cong C_3$ cyclically permutes $\chi_2,\chi_3,\chi_4$. Since $V$ is an $I$-representation, these characters occur in $V|_P$ with equal multiplicity. As $\dim V=2$, none can occur. Thus either $V|_P\cong\mathds{1}^{\oplus2}$ or $V|_P\cong\sigma$. The former cannot occur since $P$ acts faithfully on $V$. Hence $V|_P\cong\sigma$, proving the remaining assertions.
\end{proof}

Fix $E/L$ as in Lemma~\ref{lem:dyadic-elliptic}, and choose an equation $E:y^2=f(x)$ with $f\in L[x]$ a square-free cubic. For $b\in\mathcal O_{L_r}^\times$, let
\begin{equation*}
E_b:(1+\pi_Lb)y^2=f(x)
\end{equation*}
be the quadratic twist of $E$ over $L_r$ associated with $\eta_b$. For each prime $\ell\ne2$, put $V_b=V_\ell(E_b)\otimes_{\Q_\ell}\overline{\Q}_\ell$.

\begin{lemma}\label{lem:dyadic-distinct-inertia}
Each $V_b$ is wildly ramified and irreducible as an $I_L$-representation. If $\bar b\ne\bar b'$, then $V_b$ and $V_{b'}$ are non-isomorphic as $I_L$-representations.
\end{lemma}

\begin{proof}
Put $V=V_\ell(E)\otimes_{\Q_\ell}\overline{\Q}_\ell$. Then $V_b\cong V|_{G_{L_r}}\otimes\eta_b$. Lemma~\ref{lem:dyadic-elliptic} shows that $V|_{P_L}$ is irreducible, and tensoring with a quadratic character preserves irreducibility. Thus $V_b|_{P_L}$ is irreducible. The action of $P_L$ cannot be trivial, since otherwise $V_b|_{P_L}\cong\mathds{1}^{\oplus2}$ would be reducible. Thus $V_b$ is wildly ramified, proving the first assertion.

Suppose that $\bar b\ne\bar b'$ and that $V_b\cong V_{b'}$ as $I_L$-representations. Put $\eta=\eta_b\eta_{b'}^{-1}$. Then $V|_{I_L}\cong V|_{I_L}\otimes\eta$. If $h\in I_L$ acts as the identity on $V$, it acts as $\eta(h)\operatorname{id}$ on $V\otimes\eta$. Since these representations are isomorphic, $\eta(h)=1$. Thus $\eta|_{I_L}$ descends to a character $I\to\{\pm1\}$. Such a character factors through $I^{\mathrm{ab}}\cong C_3$ and is therefore trivial. This contradicts Lemma~\ref{lem:dyadic-quadratic-characters}.
\end{proof}

We now choose $b$ whose reduction generates the residue extension and apply Weil restriction.

\begin{proposition}\label{prop:dyadic-weil-realisation}
Let $r\ge1$, and choose $b\in\mathcal O_{L_r}^\times$ such that $k_{L_r}=k_L(\bar b)$. Then $A_r=\Res_{L_r/L}E_b$ has dimension $r$, and $V_\ell(A_r)$ is wildly ramified and absolutely irreducible for every prime $\ell\ne2$.
\end{proposition}

\begin{proof}
The dimension follows from $[L_r:L]=r$. Fix $\ell\ne2$, put $V=V_\ell(A_r)\otimes_{\Q_\ell}\overline{\Q}_\ell$, and let $\Phi\in G_L$ be an arithmetic Frobenius lift. For $0\le i<r$, put $b_i=\Phi^i(b)$. As in Proposition~\ref{prop:wild-weil-realisation}, Weil restriction gives
\begin{equation*}
V|_{I_L}\cong\bigoplus_{i=0}^{r-1}V_{b_i},
\end{equation*}
and $\Phi$ sends $V_{b_i}$ to $V_{b_{i+1}}$, with indices modulo $r$.

Since $k_{L_r}=k_L(\bar b)$, the residues $\bar b_i$ are distinct. Lemma~\ref{lem:dyadic-distinct-inertia} shows that these summands are irreducible and pairwise non-isomorphic under inertia. Frobenius permutes them transitively, so Lemma~\ref{lem:frobenius-transitivity} gives absolute irreducibility. Each summand is wildly ramified, and hence so is $V_\ell(A_r)$.
\end{proof}

\begin{proof}[Proof of Theorem~\ref{thm:dyadic-jacobian-realisations}\textup{(ii)}]
Choose $b\in\mathcal O_{L_g}^\times$ with $k_{L_g}=k_L(\bar b)$, and put $A_g=\Res_{L_g/L}E_b$. The elliptic curve $E_b$ has an $L_g$-rational, so Lemma~\ref{lem:compact-type-smoothing} gives a curve $C/L$ of genus $g$ with $V_\ell(\Jac(C))\cong V_\ell(A_g)$ as $G_L$-representations for every odd prime $\ell$. Proposition~\ref{prop:dyadic-weil-realisation} proves the assertion.
\end{proof}

\section{Uniform isogeny bounds over number fields}\label{sec:global}

In this section, $M$ is a number field. For each finite place $u$ of $M$, we identify $G_{M_u}$ with a decomposition subgroup of $G_M$. For an abelian variety $A/M$ and a prime $\ell$, we say that $V_\ell(A)$ is \emph{locally irreducible at $u$}, where $u\nmid\ell$, if its restriction to $G_{M_u}$ is absolutely irreducible.

To bound the degree of an isogeny $\phi:A\to B$, we compare its images on the integral Tate modules. The following lemma reduces the desired bound to an inclusion of lattices.

\begin{lemma}\label{lem:lattice-isogeny-bound}
Let $\phi:A\to B$ be an $M$-isogeny between abelian varieties of dimension $g$, and let $N\ge1$ be an integer. If $NT_\ell(B)\subseteq\phi(T_\ell(A))$ for every prime $\ell$, then $\deg\phi$ divides $N^{2g}$.
\end{lemma}

\begin{proof}
Fix a prime $\ell$ and write $N=\ell^a m$, where $a\ge0$ and $\ell\nmid m$. Since $m$ is a unit in $\Z_\ell$, we have $T_\ell(B)/NT_\ell(B)\cong(\Z/\ell^a\Z)^{2g}.$ The above inclusion gives a surjection
\begin{equation*}
T_\ell(B)/NT_\ell(B)\twoheadrightarrow
T_\ell(B)/\phi(T_\ell(A)).
\end{equation*}
The group on the left has order $\ell^{2ga}$, while the order of the group on the right is the $\ell$-part of $\deg\phi$ by \cite[Ch.~IV, Lem.~2.2(b)]{MilneAV}. Thus the $\ell$-part of $\deg\phi$ divides $\ell^{2ga}$. Since this holds for every prime $\ell$, we obtain $\deg\phi\mid N^{2g}$.
\end{proof}

\begin{lemma}\label{lem:schur-averaging}
Let $\phi:A\to B$ be an $M$-isogeny of least degree between abelian varieties of dimension $g\ge1$ over $M$. Fix a prime $\ell$ and a finite place $u\nmid\ell$ of $M$, and suppose that $V_\ell(B)$ is absolutely irreducible as a $G_{M_u}$-representation. Write $\overline G$ for the image of $G_{M_u}$ in $\mathrm{PGL}(V_\ell(B))$. If $\overline G$ is finite, then
\begin{equation*}
\frac{|\overline G|}{2g}T_\ell(B)\subseteq\phi(T_\ell(A)).
\end{equation*}
\end{lemma}

\begin{proof}
Put $L=T_\ell(B)$ and $L'=\phi(T_\ell(A))$. If $L'\subseteq\ell L$, then $A[\ell]\subseteq\ker\phi$, so $\phi$ factors through $[\ell]:A\to A$, giving an $M$-isogeny $A\to B$ of smaller degree. As a result, $L'\not\subseteq\ell L$.

Choose $x\in L'\setminus\ell L$ and extend it to a $\Z_\ell$-basis of $L$. Let $e:L\to\Z_\ell x$ be the projection along the remaining basis vectors. Then $e(L)\subseteq L'$ and $\operatorname{tr}(e)=1$.

Choose representatives in the image of $G_{M_u}$ for the elements of $\overline G$, and put $S=\sum_{h\in\overline G}heh^{-1}.$ Each representative preserves both lattices, so $S(L)\subseteq L'$. Since scalars cancel in $heh^{-1}$, each summand depends only on the image of $h$ in $\overline G$. Conjugation by any element of $G_{M_u}$ therefore permutes the summands.

After extending scalars to $\Q_\ell$, the map $S$ consequently commutes with the action of $G_{M_u}$ on $V_\ell(B)$. Absolute irreducibility and Schur's lemma give $S=c\operatorname{id}$ for some $c\in\Q_\ell$. Taking traces gives $2gc=|\overline G|$.

Finally, $cL=S(L)\subseteq L'\subseteq L$, so $c\in\Z_\ell$ and $cL\subseteq L'$, as required.
\end{proof}

We now bound $|\overline G|$ in terms of the dimension alone.

\begin{lemma}\label{lem:projective-image}
Let $A/\mathcal K$ be an abelian variety of dimension $g\ge1$ such that $V_\ell(A)$ is absolutely irreducible, where $\ell\ne p$. Then the image $\overline G$ of $G_{\mathcal K}$ in $\mathrm{PGL}(V_\ell(A))$ is finite, and its order divides $|\mathrm{GL}_{2g}(\Z/12\Z)|$.
\end{lemma}

\begin{proof}
By Lemma~\ref{lem:arithmetic-consequences}\textup{(ii)}, $A$ has potentially good reduction. Choose $m\in\{3,4\}$ prime to $p$ and put $F_m=\mathcal K(A[m])$. Then $F_m/\mathcal K$ is finite Galois and $A$ has good reduction over $F_m$ by \cite[p.~497, Cor.~2]{SerreTate}. The subgroup $I_{F_m}$ acts trivially on $V_\ell(A)$ by the N\'eron--Ogg--Shafarevich criterion, while an arithmetic Frobenius lift in $G_{F_m}$ acts as a scalar by Proposition~\ref{prop:square-residue-good-reduction}. Since these generate a dense subgroup of $G_{F_m}$, every element of $G_{F_m}$ acts as a scalar. Thus $\overline G$ is a quotient of $\Gal(F_m/\mathcal K)$.

The action on $A[m]$ gives an embedding $\Gal(F_m/\mathcal K)\hookrightarrow\mathrm{GL}_{2g}(\Z/m\Z)$. Consequently, $|\overline G|$ divides $|\mathrm{GL}_{2g}(\Z/m\Z)|$, which divides $|\mathrm{GL}_{2g}(\Z/12\Z)|$ by the Chinese remainder theorem.
\end{proof}

\begin{theorem}\label{thm:uniform-isogeny-bound!}
Let $A/M$ be an abelian variety of dimension $g\ge1$ over a number field. Suppose that, for some prime $\ell$, the representation $V_\ell(A)$ is locally irreducible at two finite places $v,w\nmid\ell$ of distinct residue characteristics.

Then there exists an integer $B_g\ge1$, depending only on $g$, such that every abelian variety $B/M$ that is $M$-isogenous to $A$ admits an $M$-isogeny $A\to B$ of degree at most $B_g$. Moreover, the $M$-isogeny class of $A$ contains at most $B_g$ $M$-isomorphism classes.
\end{theorem}

\begin{proof}
Put $N_g=|\mathrm{GL}_{2g}(\Z/12\Z)|$ and $c_g=N_g^{2g}$. Let $B/M$ be $M$-isogenous to $A$, and choose an $M$-isogeny $\phi:A\to B$ of least degree.

Let $\ell$ be any prime. Since $v$ and $w$ have distinct residue characteristics, we can choose $u\in\{v,w\}$ whose residue characteristic differs from $\ell$. By independence of the coefficient prime \cite[Cor.~1.6]{DDWeil} and invariance under $M$-isogeny, $V_\ell(B)$ is absolutely irreducible as a $G_{M_u}$-representation. Combining Lemmas~\ref{lem:projective-image} and~\ref{lem:schur-averaging} gives $N_gT_\ell(B)\subseteq\phi(T_\ell(A)).$ Since this holds for every prime $\ell$, Lemma~\ref{lem:lattice-isogeny-bound} gives $\deg\phi\mid c_g$.

Every such $B$ is $M$-isomorphic to a quotient of $A$ by a subgroup of $A[c_g]$. Thus there are at most $d_g$ possible $M$-isomorphism classes, where $d_g$ is the number of subgroups of $(\Z/c_g\Z)^{2g}$. Taking \mbox{$B_g=\max\{c_g,d_g\}$} proves both assertions.
\end{proof}

\begin{corollary}\label{cor:uniform-curves-isogeny-class}
Let $C/M$ be a curve of genus $g\ge2$ whose Jacobian satisfies the hypotheses of Theorem~\ref{thm:uniform-isogeny-bound!}. Then, up to $M$-isomorphism, the number of curves $D/M$ whose Jacobians are $M$-isogenous to $\Jac(C)$ is bounded only in terms of $g$.
\end{corollary}

\begin{proof}
Let $B/M$ be $M$-isogenous to $\Jac(C)$. Applying Lemma~\ref{lem:arithmetic-consequences}\textup{(i)} over $M_v$ gives $\End_{M_v}(B)=\Z$, and hence $\End_M(B)=\Z$. This implies that $B$ admits at most one principal polarisation defined over $M$. By Torelli's theorem over $M$ \cite[Appendix, \S1, Thms.~1 and~2]{Lauter}, $B$ is the Jacobian of at most one $M$-isomorphism class of curves of genus $g$. The result now follows from Theorem~\ref{thm:uniform-isogeny-bound!}.
\end{proof}

\end{document}

%% file: cluster_pictures.tex
\usepackage{tikz}
\usepackage{pbox}
\usetikzlibrary{fit,positioning,calc}
\def\pb#1{\pbox[c]{\textwidth}{\hfil #1\hfil}}

\definecolor{amethyst}{rgb}{0.6, 0.4, 0.8}
\definecolor{atomictangerine}{rgb}{1.0, 0.6, 0.4}
\definecolor{deeppeach}{rgb}{1.0, 0.8, 0.64}
\definecolor{eggshell}{rgb}{0.94, 0.92, 0.84}
\definecolor{lightapricot}{rgb}{0.99, 0.84, 0.69}
\definecolor{lemonchiffon}{rgb}{1.0, 0.98, 0.8}
\definecolor{roundabout}{rgb}{1.0, 0.91, 0.75}
\definecolor{atomictangerine}{rgb}{1.0, 0.6, 0.4}

\def\rootsep{0.03}               
\def\clustersep{0.06}            
\def\cnamescale{0.4}             
\def\cdepthscale{0.4}            
\def\cltopskip{1pt}              
\def\clbottomskip{1pt}           

\def\rootscale{0.5}   \def\rootcolor{amethyst}
\def\rootscaleA{0.7}  \def\rootcolorA{yellow}
\def\rootscaleB{0.5}  \def\rootcolorB{green}
\def\rootscaleC{0.5}  \def\rootcolorC{red}
\def\rootscaleD{0.5}  \def\rootcolorD{blue}
\tikzset{
  clA/.style = {very thick,black},
  clB/.style = {thick,purple}
}

\tikzset{
  root/.style = {circle,scale=\rootscale,ball color=\rootcolor},
    rc/.style 2 args = {right=#1*1.5*\clustersep of {#2.east|-first},root}, rr/.style = {right=\rootsep of {#1.east|-first},root},
  roott/.style = {circle,inner sep=-2pt,minimum size=5pt,black,font=\ttfamily\footnotesize},
    rct/.style 2 args = {right=#1*1.5*\clustersep of {#2.east|-first},roott}, rrt/.style = {right=\rootsep of {#1.east|-first},roott},
  rootA/.style = {circle,scale=\rootscaleA,ball color=\rootcolorA},
    rcA/.style 2 args = {right=#1*1.5*\clustersep of {#2.east|-first},rootA}, rrA/.style = {right=\rootsep of {#1.east|-first},rootA},
  rootB/.style = {circle,scale=\rootscaleB,ball color=\rootcolorB},
    rcB/.style 2 args = {right=#1*1.5*\clustersep of {#2.east|-first},rootB}, rrB/.style = {right=\rootsep of {#1.east|-first},rootB},
  rootC/.style = {circle,scale=\rootscaleC,ball color=\rootcolorC},
    rcC/.style 2 args = {right=#1*1.5*\clustersep of {#2.east|-first},rootC}, rrC/.style = {right=\rootsep of {#1.east|-first},rootC},
  rootD/.style = {circle,scale=\rootscaleD,ball color=\rootcolorD},
    rcD/.style 2 args = {right=#1*1.5*\clustersep of {#2.east|-first},rootD}, rrD/.style = {right=\rootsep of {#1.east|-first},rootD},
  cluster/.style = {draw=blue!70,thick,rounded corners,inner sep=22*\clustersep,outer xsep=22*\clustersep,fit=#1},
  clabel/.style  = {anchor=west,scale=\cdepthscale,black,inner sep=0,outer xsep=1,outer ysep=0},
  clabelL/.style = {above right=-\clustersep of #1t.north east,clabel},
  clabelD/.style = {below right=-\clustersep of #1t.south east,clabel},
  clouter/.style = {inner sep=0,outer sep=0,fit=#1}
}

\def\Cluster #1 = #2;{\node[cluster=#2] (#1) {};}
\def\ClusterL #1[#2] = #3;{
  \node[cluster=#3] (#1t) {}; \node[clabelL=#1] (#1l) {$#2$}; \node[clouter=(#1t)(#1l)] (#1) {};}
\def\ClusterD #1[#2] = #3;{
  \node[cluster=#3] (#1t) {}; \node[clabelD=#1] (#1d) {$#2$}; \node[clouter=(#1t)(#1d)] (#1) {};}
\def\ClusterLD #1[#2][#3] = #4;{
  \node[cluster=#4] (#1t) {}; \node[clabelL=#1] (#1l) {$#2$}; 
  \node[clabelD=#1] (#1d) {$#3$}; \node[clouter=(#1t)(#1l)(#1d)] (#1) {};}
\def\ClusterLDName #1[#2][#3][#4] = #5;{
  \node[cluster=#5] (#1t) {}; \node[clabelL=#1] (#1l) {$#2$}; 
  \node[clabelD=#1] (#1d) {$#3$}; 
  \node[scale=\cnamescale,above=\clustersep/3 of #1t,inner sep=0, outer sep=0] (#1n) {$#4$}; 
  \node[clouter=(#1l)(#1d)(#1t)] (#1) {};}

\newcommand{\Root}[4][]{
  \ifx\relax#2\relax\node[rr#1=#3] (#4) {};\else\node[rc#1={#2}{#3}] (#4) {};\fi}
\newcommand{\RootT}[5][]{
  \ifx\relax#2\relax\node[rrt#1=#3] (#4) {#5};\else\node[rct#1={#2}{#3}] (#4) {#5};\fi}

\def\frob(#1)(#2){\path[draw,thick,shorten <=-22*\clustersep,shorten >=-22*\clustersep](#1.east)--(#2.west|-#1){};}

\long\def\clusterpicture#1\endclusterpicture{\pb{\vbox to \cltopskip{\vfill}\\%
  \begin{tikzpicture}\node[coordinate] (first) {};#1\end{tikzpicture}\\[-11pt]\vbox to \clbottomskip{\vfill}}}   
\long\def\clusterpictureopt#1#2\endclusterpicture{\pb{\vbox to \cltopskip{\vfill}\\%
  \begin{tikzpicture}[#1]\node[coordinate] (first) {};#2\end{tikzpicture}\\[-11pt]\vbox to \clbottomskip{\vfill}}}